\documentclass[12pt]{article}

\usepackage{amssymb, amsmath}

\newtheorem{Thm}{Theorem}
\newtheorem{Lem}[Thm]{Lemma}
\newtheorem{Cor}[Thm]{Corollary}

\newtheorem{Rmk}{Remark}
\def\qed{\kern 6pt\hbox{\vrule\vbox to 6pt{\hrule width
                            6pt\vfil\hrule}\vrule}}

\newcommand{\R}{\mathbb{R}}

\def\cal{\mathcal}

 \def\exp{\mathop{\rm
exp}\nolimits}

\title{On the volume product of planar polar convex bodies --- 
lower estimates with stability\\
{\rm{Stud. Sci. Math. Hungar. 50 (2) (2013), 159-198\\
DOI:10.1556.SScMath.49.2012.2.12.197}}}

\author{
K. J. B\"or\"oczky* \\
Alfr\'ed R\'enyi Mathematical Institute, \\
Hungarian Academy of Sciences, \\
H-1364 Budapest, Pf. 127, HUNGARY, and  \\
Universitat Polit\'ecnica de Catalunya, Barcelona Tech., SPAIN  \\
boroczky.karoly@renyi.mta.hu  \\
http://www.renyi.hu/$\,\,\widetilde{}\,$carlos  \\
E. Makai, Jr.**  \\
Alfr\'ed R\'enyi Mathematical Institute, \\ 
Hungarian Academy of Sciences,  \\
H-1364 Budapest, Pf. 127, HUNGARY  \\
makai.endre@renyi.mta.hu  \\
http://www.renyi.hu/$\,\,\widetilde{}\,$makai  \\
M. Meyer***  \\
\'Equipe d'Analyse et de Math\'ematiques Appliqu\'ees, \\
Universit\'e de Marne-la-Vall\'ee,
Cit\'e Descartes, \\
5, boulevard Descartes, Champs-sur-Marne, \\
77454 Marne-la-Vall\'ee Cedex 2, FRANCE \\
mathieu.meyer@univ-mlv.fr \\
http://umr-math.univ-mlv.fr \\
S. Reisner****  \\
Department of Mathematics, University of Haifa, \\
Haifa 31905, ISRAEL \\
reisner@math.haifa.ac.il \\
http://math.haifa.ac.il/reisner}

\begin{document}

\thanks{Research (partially) supported by 
Hungarian National Foundation for Scientific Research,  
grant nos. K68398, K75016, FP7 IEF grant GEOSUMSETS. 
$^{**}$Research (partially) supported by 
Hungarian National Foundation for Scientific Research, 
grant nos. K68398, K75016, K81146. 
$^{***}$The research has been partially supported by 
the France-Israel Research Network Program in Mathematics. 
$^{****}$The research has been
partially supported by 
the France-Israel Research Network Program in Mathematics, and 
Discrete and Convex Geometry, 
Marie Curie Host Fellowships for Transfer of Knowledge, 
MTKD-CT-2004-014333. 
Part of the work was done when S. Reisner visited R\'enyi Institute 
in the summer of 2008. He is indebted to the Institute 
for the hospitality and the excellent research conditions.}


\maketitle

\begin{abstract} Let $K \subset {\mathbb R}^2$ 
be an $o$-symmetric convex body, and $K^*$ its polar body. 
Then we have $|K|\cdot |K^*|
\ge 8$, with equality if and only if $K$ is a parallelogram.
($| \cdot |$ denotes volume).
If $K \subset {\mathbb R}^2$ 
is a convex body, with $o \in {\text{int}}\,K$, then
$|K|\cdot |K^*|
\ge 27/4$, with equality if and only if $K$ is a triangle and $o$ is its
centroid.
If $K \subset {\mathbb R}^2$ 
is a convex body, then we have
$|K| \cdot |[(K-K)/2)]^* | \ge 6$, with equality if and only if $K$ is a
triangle. 
These theorems are due to Mahler and Reisner, Mahler and Meyer, 
and to Eggleston, respectively. We show an analogous theorem:
if $K$ has $n$-fold rotational symmetry about $o$, then $|K|\cdot |K^*|
\ge n^2 \sin ^2 ( \pi /n)$, 
with equality if and only if $K$ is a regular $n$-gon of centre $o$.  
We will also give stability variants of these four inequalities, both for the
body, and for the centre of polarity. 
For this 
we use the Banach-Mazur distance
(from parallelograms, or triangles), or its analogue with similar copies
rather than affine transforms (from regular $n$-gons),
respectively. The stability variants are sharp, up to constant factors.
We extend the inequality $|K|\cdot |K^*| \ge n^2 \sin ^2 ( \pi /n)$
to bodies with $o \in {\text{int}}\,K$, 
which contain, and are contained in,
two regular $n$-gons, the vertices of the contained $n$-gon
being incident to the sides of the containing $n$-gon. 
Our key lemma is a stability estimate for the area product of two
sectors of convex bodies polar to each other.
To several of our
statements we give several proofs; in particular, we give a new proof for the
theorem of Mahler-Reisner. 

\end{abstract}



RUNNING AUTHORS: K. J. B\"or\"oczky, E. Makai, Jr., M. Meyer, S. Reisner

RUNNING TITLE: Volume product of planar convex bodies


\section{Notations}

We write
$o$ for the origin,
$\langle \cdot,\cdot \rangle$ for the scalar product,
$\|\cdot\|$ for the Euclidean norm, $[x_1,\ldots,x_k]$
for the convex hull of $ \{ x_1,\ldots,x_k \} $, and 
$| \cdot |$ for the volume (area in ${\mathbb R}^2$).
We write vert\,$( \cdot )$, int\,$( \cdot )$, $\partial ( \cdot )$, for the
vertex set of a convex polytope, or interior, or boundary of a set in
${\R}^d$. 

A convex body in $\R^d$ is a compact convex set with non-empty interior.
If $o\in {\rm int}\,K$, then
its polar (w.r.t. the unit sphere with centre $o$) is
$$
K^*=\{x\in\R^d:\, \forall y \in K\,\,\,\, \langle x,y\rangle\leq 1\}.
$$
If $A:\R^d\to\R^d$ is a non-singular linear map,
then $(AK)^*=(A^{-1})^*K^*$, where $(A^{-1})^*$ is the transpose
of the inverse of $A$.
It is known (Santal\'o \cite{San}, or Meyer-Pajor \cite{MP90}),
that there exists a unique point
$s(K) \in{\rm int}\,K$, called {\it{Santal\'o point of $K$}}, such that
$$
|\left( K-s(K) \right) ^*|=\min \, \{ |(K-z)^*| : z \in {\rm int}\,K \}.
$$
Additionally, the origin is the centroid of $\left( K-s(K) \right) ^*$.
The uniqueness and the affine invariance of the Santal\'o point yields that
$s(K)=o$ if $K$ is $o$-symmetric, or if $d=2$ and $K$
has $n$-fold symmetry about $o$ for some $n \geq 3$.

For convex bodies $K,L \subset {{\mathbb R}^d}$, 
the {\it{Banach-Mazur distance}}
$\delta _{BM}(K,L)$ is $\min \, \{ \lambda _2 / \lambda _1 \mid \lambda _1,
\lambda _2 \in (0, \infty ),\,\, \exists 
{\text{ affinity }} A, \,\, \exists x \in {{\mathbb R}^d}, \,\, \lambda _1 A(K)
\subset L \subset \lambda _2 A(K) +x \} $. If we allow for $A$ only
similarities, then we obtain the definition of $\delta _{BM}^s
(K,L)$. (Clearly, $\delta _{BM}(K,L) \le \delta _{BM}^s (K,L)$. If 
both $K,L$ are $o$-symmetric, or $d=2$ and both
have $n$-fold rotational symmetry about $o$,
with $n \ge 3$ an integer,
then in the definition of $\delta _{BM}(K,L)$, or $\delta _{BM}^s (K,L)$, we
may assume $x=o$.) We will write $T,P,R_n$ for a triangle, parallelogram, or
regular $n$-gon, respectively.

We write $\kappa _d=\pi ^{d/2}/\Gamma (d/2+1)$ for the volume of the
unit ball in the $d$-dimensional Euclidean space. 


\section{Introduction}

{\bf{Generalities.}}
Let $K \subset {\mathbb R}^d$ 
be a convex body, with $o \in {\text{int}}\,K$. 
Blaschke \cite {Bl} was the 
first who considered the so called {\it{volume product}} $|K| \cdot
|K^*|$ of the body $K$, and proved that for $d \le 3$, and $o$ the barycentre
of $K$, its maximum is attained,
e.g., if $K$ is an ellipsoid. He was motivated by the investigation of
the affine geometry of
convex bodies, e.g., of the so called {\it{affine surface area}} 
(a definition
cf. in \cite{L}, or \cite{Bo}), that is intimately related to the volume
product (cf. \cite{L}, \cite{Bo}).
The volume product is invariant under non-singular
linear transformations, cf. \cite{L},
p. 109. The investigation of the question of the lower estimate
of the volume product was initiated by
Mahler \cite{Mah38}, \cite{Mah39}. He had in view applications in
the geometry of numbers (i.e., investigation of the relation of convex, or
more generally, of star-bodies, 
and lattices, i.e., non-singular linear images
of ${\mathbb Z}^d$ in ${\mathbb R}^d$). The volume product, in particular, 
for
$o$-symmetric $K$, is a basic
quantity, that later has arisen in several branches of mathematics,
cf. later in this introduction.

{\bf{Lower bound, $o$-symmetric case.}}
For a while we suppose that $K$ is $o$-symmetric.
Mahler \cite{Mah39}, for $d \ge 2$, conjectured for the volume product
the lower bound $4^d/d!=d^{-d}\left( 4e +o(1) \right) ^d$, 
and proved the lower bound $4^d/(d!)^2$. 
It is usually credited to Saint Raymond \cite{SR} 
that this conjectured lower bound is attained
not only for parallelotopes and cross-polytopes. However, this had
already been observed by Guggenheimer \cite{Gu} 
some years earlier, where the way of obtaining all examples of
\cite{SR} had already been 
described. These examples are the following. Beginning with $[-1,1] \subset
{\mathbb R}$, we define inductively convex bodies in ${\mathbb R}^d$, 
from examples in lower dimensions: if $d=d_1+d_2$ is an arbitrary
decomposition of $d$ as a sum of positive integers $d_1,d_2$, then
for the already defined bodies in ${\mathbb R}^{d_i}$ we take either their
Minkowski sum, or the convex hull of their union. The Banach spaces with
these unit balls are called {\it{Hansen-Lima spaces}}, and their unit balls are
called {\it{Hansen-Lima bodies}}.
They are called also {\it{Hanner bodies}}. 
Since they were introduced already by
Hammer \cite{Ha} in 1956, while Hansen-Lima \cite{HL} dates to 1981, it would
be more correct to call them {\it{Hanner-Hansen-Lima spaces, and bodies}}.
Also the well known book of Gr\"unbaum on polytopes \cite{Gru} 
calls these bodies {\emph{Hanner
polytopes}}. Moreover,
\cite{Ha} and \cite{HL} proved two characterizations of these convex bodies
$K$: (1) (a Helly type property): $K$ is $o$-symmetric, and if three translates
of $K$ pairwise intersect, then their intersection is not empty. (2): $K$ is
an $o$-symmetric convex polytope, and for any two disjoint faces (of any
dimensions $ \le d-1$), say, $K_1,K_2$, of $K$, there are two distinct parallel
supporting hyperplanes $\Pi _1$ and $\Pi _2$ of $K$, such that $K_1 \subset
\Pi _1$, and $K_2 \subset \Pi _2$.
\cite{SR} conjectured that the
volume product attains its minimum exactly for the Hanner-Hansen-Lima bodies. 
(However, the claim of \cite{Gu} that its author
settled the $3$-dimensional case is incorrect.) 

Mahler \cite{Mah38}
proved {\it{the sharp lower bound $|K| \cdot |K^*| \ge 8$ for $d=2$ and $K$
$o$-symmetric.}} 
Morerover, \cite{Mah38} showed that, {\it{for $K$ a
polygon, the lower bound is attained, for the $o$-symmetric case, 
if and only if $K$ is a parallelogram}}.

The above lower estimate of \cite{Mah39} for
${\mathbb R}^d$, for the $o$-symmetric case,
was sharpened to $2^d \kappa _d /(d!d^d)^{1/2}$ by Dvoretzky-Rogers 
\cite{DR}, and to 
$\kappa _d^2/d^{d/2}$ by Bambah \cite{Bam}. Then it has become clear that the
volume product is very important in functional analysis, where it is just the
product of the volumes of the unit balls of a finite dimensional Banach space
and its dual. This has importance in the so called local theory of Banach
spaces, i.e., the asymptotic 
study of finite dimensional Banach spaces, of high dimension, cf. Pisier's
book \cite{Pi}. A number of
other geometric characteristics of these Banach spaces have a 
connection to the
volume product. Therefore functional analysts became strongly interested in
the subject, which resulted in ever better {\it{lower estimates}}, namely  
$d^{-d} (\log d)^{-d} \cdot {\text{const}}^d 
$ by Gordon-Reisner \cite{GR} and later by Kuperberg, G. \cite{K92},
and to $d^{-d} \cdot {\text{const}} ^d  $ 
by Bourgain-Milman \cite{BoMi}
(with an unspecified constant). Quite recently
$(\pi /4)^{d-1} 4^d/d!
=d^{-d}\left( e \pi +o(1) \right) ^d$ {\it{was proved 
by Kuperberg, G.}} \cite{K08}. Observe that the quotient of
Kuperberg G.'s estimate and the
conjectured minimum is $(\pi / 4 + o(1))^d$. The paper Nazarov \cite{N} 
proved the bit weaker estimate $(4^d/d!)( \pi /4)^{3d}$, using an unexpected
connection of the volume product problem to the theory of functions of
several complex variables.

{\it{A class of $o$-symmetric convex bodies in ${\mathbb R}^d$, 
for which the lower bound $4^d/d!$ is
known, is the class of (non-singular) linear images of convex bodies
symmetric with respect
to all coordinate hyperplanes}} (also called {\it{unconditional convex 
bodies}}), 
cf. Saint Raymond \cite{SR}.
The equality cases were proved
by Meyer \cite{Me86} and Reisner \cite{R87} --- these are just the
Hanner-Hansen-Lima bodies. The combinatorial aspect of the proof was later
studied by Bollob\'as-Reader-Radcliffe \cite{BLR}. 
Actually \cite{SR} proved this inequality for a larger
class of $o$-symmetric convex bodies. These are the ones,
for which the associated norm satisfies the following. There exists a
base, such that for the coordinates $x_1, \ldots ,x_d$ w.r.t. this base,
the projections $(x_1, \ldots ,x_d) \to (x_1, \ldots x_{i-1},x_{i+1}, \ldots
x_d)$, where $1 \le i \le d$, are contractions. 
Moreover, \cite{SR} also extended his inequality, for unconditional convex
bodies, in the following way. Let $k \ge 2$ be an integer, let
an unconditional norm $\| \cdot \| $
on ${\mathbb R}^k$ be given (i.e., the unit ball is unconditional), and let
$d_1,...,d_k \ge 1$ be integers. Let
$K_i \subset {\mathbb R}^{d_i}$ be $o$-symmetric convex bodies, which are the
unit balls of norms $\| \cdot \| _i$. We consider $\prod _{i=1}^k {\mathbb R}
^{d_i}$, with the norm $\| ( \| x_i \| _i ) \| $ (where $x_i \in {\R}^{d_i}$), 
where we consider $\| \cdot
\| _i $ as fixed, and $\| \cdot \| $ as variable. Then the volume product of
the unit ball of this norm attains its minimum, e.g., for the cases, when 
$\| ( \lambda _1, ..., \lambda _k) \| $ 
equals $\sum _i |\lambda _i|$, or $\max _i |\lambda _i|$.

{\it{Mahler's conjecture in the $o$-symmetric case, 
together with the conjecture about the equality cases, 
is also proved for convex polytopes
with (at most) $2d+2$ vertices or facets, for}} $d \le 8$, cf. Lopez-Reisner 
\cite{LR}.

{\it{Mahler's conjecture is also proved for 
zonoids $K$ in}} ${\mathbb R}^d$ (i.e., {\it{limits in
the Hausdorff-metric of finite sums of segments}}), {\it{with centre at $o$,
and with}} ${\text{int}}\,K \ne \emptyset $. This is due to 
Reisner \cite{R85}, \cite{R86}, and in these papers it is also proved that
{\it{the lower bound for zonoids is attained if and only if $K$ is a
parallelotope, with centre at}} $o$. 
Later, a simpler proof was given by
Gordon-Meyer-Reisner \cite{GMR}. 
Observe that this settles the case of equality for
$o$-symmetric convex bodies in ${\mathbb R}^2$, 
since each such body is a zonoid. Both \cite{R85}, \cite{R86} use the
connection of the volume product problem with stochastic geometry (geometric
probability), as is done
also later in B\"or\"oczky, K. J.-Hug \cite{BH}, in another context. 
A variety
of other connections to geometric probability are contained in Thompson's 
book
\cite{Th96}, in particular in Ch. VI.
\cite{R85} also gave an analogue of the last mentioned
Saint Raymond's theorem: if each $K_i$, there considered,
is either a zonoid, or the polar of a
zonoid, then $|K| \cdot |K^*| \ge 4^d/d!$. Reisner 
\cite{R87} clarified the equality cases in the last mentioned Saint Raymond's
theorem (\cite{SR}): this is the case if and only if 
$\| \cdot \|$ is a norm of a Hanner-Hansen-Lima space.

Barthe-Fradelizi in the paper \cite{BF}
proved that {\it{if $K$
is a convex body and $P$ is a regular convex polytope in $\R^d$ such that
the origin is the centroid of $P$, and $K$ has all the symmetries
of $P$}} --- thus the origin is also their common Santal\'o point ---
{\it{then $|K|\cdot|K^*|\geq |P|\cdot|P^*|$, with equality if and only if 
$K$ is a dilate of $P$ or of $P^*$}}. Its particular case for $d=2$ 
is given in a bit stronger form in our 
Corollary~\ref{nfold}; however, we have as well a stability variant of our
Corollary~\ref{nfold}, in our Theorem~\ref{nfoldstab}.
\cite{BF} has some generalizations of the above cited inequality as well.
Let $d=d_1+...+d_k$ be a decomposition of $d$ into positive integers. Let us
have in ${\mathbb R}^d$ a convex body $K$, whose symmetry group ${\cal
O}(K)$ (i.e., the group of congruences mapping $K$ onto $K$) contains
the group ${\cal O}(P_1) \times ... \times {\cal O}(P_k)$, where each $P_i$ is
either a regular convex polytope or a ball in ${\mathbb R}^{d_i}$, of centre
$o$. Then we have $V(K)V(K^*) \ge V(P_1 \times ... \times P_k)
V \left( (P_1 \times ... \times P_k)^* \right) $. Here the equality
cases are not clarified. However, there are a lot of equality cases: for each
$i$, we may take $P_i$ or $P_i^*$, and may construct
from these inductively new bodies, like in case of the 
Hanner-Hansen-Lima bodies.
We remark that our \S 4 has a considerable overlap with \cite{BF}.

{\bf{Upper bound, $o$-symmetric case.}}
{\it{In the $o$-symmetric case, 
the sharp upper bound is $\kappa _d ^2=d^{-d}\left( 2e \pi +o(1) \right) ^d$,
and it
is attained if and only if $K$ is an $o$-symmetric ellipsoid}},
which is due to Blaschke \cite{Bl} ($d \le 3$) and Santal\'o \cite{San} 
(for general $d$), 
{\it{with the equality case proved by Saint Raymond}} \cite{SR}. 
Ball \cite{Bal} and Meyer-Pajor \cite{MP89} pointed
out that a proof of the inequality can be given by Steiner symmetrization:
namely that {\it{Steiner symmetrization does not decrease}} $|K| \cdot |K^*|$.
A number of simplifications of these
proofs has appeared.

Fradelizi-Meyer 
\cite{FM07} also considers the upper estimate for the volume product for
measures other than the Lebesgue measure.

{\bf{Lower bound, general case.}}
As an application to the original volume product problem, Meyer-Reisner
\cite{MR98} gives
the following statement. 
{\it{If all non-empty intersections of $K$ with horizontal
hyperplanes are positive
homothets of a given $(d-1)$-dimensional convex body $L$, and
these intersections have their 
Santal\'o points}} (taken in their affine hull)
{\it{on a line,
then $|K| \cdot |K-s(K)|/(|L| \cdot |L-s(L)|)$ attains its minimum}}
$(d+1)^{d+1}/d\,^{d+2}$ (that is independent of $L$), {\it{if and
only if $K$ is a cone, with base a translate of}} $L$. ({\it{Examples}} 
of such bodies
{\it{are bodies rotationally symmetric about the $x_d$-axis}}.) 

{\it{Mahler}} \cite{Mah38} {\it{proved the sharp lower 
bound \,$|K| \cdot |K^*| \ge 27/4$,
for $K \subset {\mathbb
R}^2$ a convex body, with}} $o \in {\text{int}}\,K$.
Morerover, \cite{Mah38} showed that, {\it{for $K$ a
polygon, the lower bound is attained, for the
case}} $o \in {\text{int}}\,K$,
{\it{if and only if $K$ is 
a triangle with barycentre at $o$}}. Later Meyer \cite{Me91} 
showed that {\it{for}} the case $d=2$ {\it{and}} $o \in {\text{int}}\,K$,
{\it{the lower bound is attained only for triangles, with barycentre 
at}}\, $o$.
A simpler proof of this is contained in Meyer-Reisner \cite{MR06}, 
Theorem 15.

For the case $o \in {\text{int}}\,K$, Mahler \cite{Mah39} 
conjectured that $|K| \cdot |K^*| 
\ge (d+1)^{d+1}/(d!)^2 \sim {\text{const}} \cdot d^{-d}e^{2d}$, 
where equality stands only for a simplex with barycentre at
$o$. {\it{The lower bound}} $(d+1)^{d+1}/\left( d^d(d!)^2 \right) $ 
is due to Mahler \cite{Mah50}, 
that was sharpened
to $\kappa _d^2/(d!)^2$ by Bambah \cite{Bam}, 
to ${\text{const}} ^d \cdot d^{-d}$ 
by Bourgain-Milman \cite{BoMi}
(with an unspecified constant), while 
$\left( \pi /(2e) \right) ^{d-1} (d+1)^{d+1}/(d!)^2$ 
{\it{has been recently proved by
Kuperberg, G.}} \cite{K08}. Observe that the quotient of
this estimate and the
conjectured minimum is $\left( \pi / (2e) + o(1) \right) ^d$.
{\it{Mahler's conjecture, for the asymmetric case, together with the conjecture
about the equality cases,
is proved for convex polytopes with at most $d+3$ vertices or facets}}, 
cf. Meyer-Reisner \cite{MR06}, Theorem 10.

Observe that, for bodies $K$ having all symmetries of a regular simplex, the
Barthe-Fradelizi result, from \cite{BF}, above cited,
implies Mahler's conjecture, in the asymmetric case,
together with the cases of equality. 
\cite{BF} contains also the following. {\it{Let us have in ${\R}^d$
a convex body $K$, such that the group of affinities preserving $K$
contains not necessarily orthogonal symmetries w.r.t. affine hyperplanes
$H_1,...,H_m$, where $\cap _{i=1}^m H_i$ is a one-point set. Then we have
$|K| \cdot | \left( K-s(K) \right) ^* | \ge (d+1)^{d+1}/(d!)^2$, with equality 
if and only if $K$ is a simplex.}} In other words, for these bodies Mahler's
conjecture, in the asymmetric case, 
about the lower bound of the volume product is true.
 
Cf. also, e.g., the recent papers Hug \cite{H}, Klartag-Milman \cite{KM},
Campi-Gronchi \cite{CG},
Meyer-Reisner \cite{MR06}, Fradelizi-Meyer \cite{FM07}, 
Fradelizi-Meyer \cite{FM08a}, Fradelizi-Meyer \cite{FM08b}, B\"or\"oczky, K. J.
\cite{Bo}, Lin, Youjiang and Leng, Gangsong \cite{LL}, 
B\"or\"oczky, K. J.-Hug \cite{BH}, 
Fradelizi-Gordon-Meyer-Reisner \cite{FGMR}, Fradelizi-Meyer
\cite{FM}, and the references therein.

{\bf{Upper bound, general case.}}
{\it{One has for 
$|K| \cdot | [ K-s(K) ] ^* | $ the
upper estimate $\kappa _d ^2$, with equality if and only if $K$ is an
ellipsoid}}, cf. Blaschke \cite{Bl}, Santal\'o \cite{San} 
for the inequality, and Petty \cite{Pe}, Meyer-Pajor \cite{MP90} 
for the cases of
equality. Again, \cite{MP90} used for the proof, among others, 
Steiner's symmetrization, but in a more involved manner, than in the
$o$-symmetric case. Recently Artstein-Avidan-Klartag-Milman \cite{A-AKM} and
Meyer-Reisner \cite{MR06} showed that Steiner
symmetrization proves the Blaschke-Santal\'o inequality, namely that 
{\it{Steiner
symmetrization does not decrease}} $|K| \cdot |(K-s(K))^*|$ (the case when $K$
is $o$-symmetric, was cited above). \cite{MR06}
proved in this way also the case of equality.
Actually {\it{the same upper estimate $\kappa _d ^2$ holds for 
$|K| \cdot | [K-b(K)] ^* |$, 
where $b(K)$ is the barycentre of $K$,
and again with equality if and only if $K$ is an
ellipsoid}}, cf. \cite{L}, p. 165. Actually, if $s(K)$, or $b(K)$, is $o$,
then $b(K^*)$, or $s(K^*)$, is $o$, respectively, cf. \cite{L}, p. 165, 
which explains the symmetric
role of the Santal\'o point, and the barycentre.

A general reference to these problems, and their connections to other affine
inequalities for convex bodies, 
is Lutwak \cite{L}. A more recent survey on the volume product
is Thompson \cite{Th07}.

{\bf{Eggleston-Zhang type problems.}}
For another generalization 
of the volume product, from the $o$-symmetric case to the general case,
Eggleston \cite{E} proved the following. {\it{If
$K \subset {\mathbb R}^2$ is a convex body, then $|K| \cdot 
| [(K-K)/2] ^* |
\ge 6$, with equality if and only if $K$ is a triangle}}. 

A generalization of this to ${\mathbb R}^d$, however
not for polar bodies, but for polars of
projection bodies, was given by Zhang \cite{Z}: 
his inequality is $|K|^{d-1} \cdot | 
(\Pi K)^*) | \ge \binom{2d}{d}d^{-d}$, {\it{with equality if and only if $K$
is a simplex}}. (The {\it{projection body $\Pi K$ 
of a convex body}} $K \subset {\mathbb R}^d$ is the 
$o$-symmetric convex body --- actually a zonoid --- whose support function at a
point $u \in S^{d-1}$ is given as the $(d-1)$-volume of the orthogonal
projection of $K$ to the linear subspace orthogonal to $u$. Observe
that for $d=2$ the bodies $\Pi K$ and $K-K$ can be obtained from each other by
a rotation through $ \pi /2$ about the origin, hence their polars 
have equal areas.)
B\"or\"oczky, K. J. \cite{Bo05}, Theorem 3 
proved an almost sharp stability version of
this inequality: {\it{for $S$ a simplex, $|K|^{d-1} \cdot 
| (\Pi K)^*) | \le \binom{2d}{d}d^{-d}(1+\varepsilon)$ implies}} $\delta
_{BM} (K, S) \le 1+ {\text{const}}_d \cdot \varepsilon ^{1/d}$, {\it{while
the actual error term cannot be less than}} 
${\text{const}}_d \cdot \varepsilon ^{1/(d-1)}$ (\cite{Bo05}, Example 19),
which quantity is conjectured
to be the exact order of the error term.

For the original question about the lower estimate of 
$|K| \cdot | [(K-K)/2] ^* |$, 
for $K \subset {\mathbb R}^d$ a convex
body, the sharp lower bound is conjectured to be $(d+1)2^d/ d! \sim 2^d e^d
d^{-d} \left( 1+o(1) \right) ^d$, with equality
for $K$ a simplex, cf. \cite{Mak78}. 
(A calculation, that for $K$ a simplex we have
equality, is given 
in \cite{MM}.) This quantity occurs in a number of problems of 
the theory of packings and coverings, and more generally in density 
estimates of sytems of convex sets
(for the non-symmetric case seemingly even more than the
original volume product), 
cf. e.g., \cite{Mak78}, \cite{Mak87}, \cite{MM96} Theorem 5.2, Remark 5.3.
Since $|K| \cdot | [(K-K)/2] ^* |=
[|K|/| (K-K)/2 | ] \cdot [ | (K-K)/2 | \cdot 
| [(K-K)/2] ^* | ] $,
Kuperberg G.'s result and the {\it{difference body inequality}} 
(i.e., $|K-K|/|K| \le {\binom{2d}{d}}$, cf. Rogers-Shephard \cite{RS}) imply 
$|K| \cdot |  [ (K-K)/2 ] ^* | 
\ge 
( \pi /4)^{d-1} 8^d/ \left( d!{\binom{2d}{d}} \right) 
\sim 
d^{-d} (e \pi /2)^d  \left( 1+o(1) \right) ^d$.
Observe that the quotient of this value and the conjectured value is 
$\left( \pi /4 +o(1) \right) ^d$.

{\bf{Florian's inequalities.}} 
A question of another character was treated by Florian in \cite{F96} 
and \cite{F98}. 
He investigated {\it{convex bodies in ${\mathbb R}^2$, 
contained in the unit circle about $o$}}, and showed the {\it{sharp estimate
$|K|+|K^*| \ge 6$, attained for a square inscribed to the unit circle}}. 
He gave as well a stability result in a more special
case. See references to earlier results of this type as well 
in \cite{F96} and \cite{F98}.

{\bf{Local and global stability results.}}
{\it{A stability version of the
Blaschke-Santal\'o inequality, for $d \ge 3$, is proved by B\"or\"oczky, 
K. J.}}
\cite{Bo} (stability meant {\it{for the Banach-Mazur distance}}).
{\it{For $d=2$ the same is done in Ball-B\"or\"oczky, K. J.}} \cite{BB}, 
{\it{B\"or\"oczky, K. J.-Makai, Jr.}} \cite{BoMa}. 

After essentially finishing our paper we were informed from the paper
Nazarov-Petrov-Ryabogin-Zvavitch \cite{NPRZ}
about the following theorem. {\it{For $d \ge 2$ an integer there exist
$\varepsilon _d >0$ and $c_d>0$ with the following properties.
If the Banach-Mazur distance of an $o$-symmetric convex body $K \subset
{\mathbb{R}}^d$ from the class of
parallelotopes is $1+ \varepsilon  \in (1, 1+ \varepsilon _d]$, 
then the volume product
$|K| \cdot |K^*|$ is at least $[4^d/d!](1+c_d \varepsilon )$.
Here the order of the error term is optimal.}} Together with the paper
B\"or\"oczky, K. J.-Hug
\cite{BH} (which calls the attention to the fact that,
although \cite{NPRZ} states its theorem in the form that parallelotopes are
strict local minima, the proof in
\cite{NPRZ} actually gives this stronger, namely, 
stability variant, cited above; cf. \cite{NPRZ}, \S 4),
this gives the following. 
{\it{For the case of $o$-symmetric zonoids $K$
in $\mathbb{R}^d$, with}} ${\text{int}}\,K
\ne \emptyset $, {\it{in particular,
for $o$-symmetric convex bodies in $\mathbb{R}^2$,
we have global stability of the parallelotopes}}, i.e.,
the above inequality, without a restriction of the form
$0 < \varepsilon \le \varepsilon _d$. For $\mathbb{R}^2$, 
this is our Theorem~\ref{MahlerReisnerstab}, without the
specification of the coefficient of $\varepsilon $ in the lower estimate.
Once more, {\it{the order of the error term is optimal}}. 

Since optimality of the order of the
above two error terms was not proved in \cite{NPRZ} or \cite{BH}, 
we show it. Of course, it suffices to deal with the case of $o$-symmetric
zonoids only, for which we give the following example.
For $d=2$ we take $[-1,1]^2$, and cut off small isosceles right triangles of
legs $\varepsilon $ at each vertex. 
For $d \ge 3$ we take the product of this example with
$[-1,1]^{d-2}$. Thus we obtain an $o$-symmetric zonoid, $K$, say. 
Then $| K | \cdot | K^* | =
(4^d/d!)\left( 1+ c_1 \varepsilon +O(\varepsilon ^2) \right) $, 
for some $c_1>0$. 
Clearly $\delta _{BM}(K,[-1,1]^d) \le 
1+ c_2 \varepsilon +O(\varepsilon ^2)$, for some $c_2>0$. 
Now we estimate $\delta _{BM}(K,[-1,1]^d) =
\delta _{BM} \left( K^*, {\text{conv}} \, \{ \pm e_i 
\mid 1 \le i \le d \} \right) $
from below, by $1+ c_3 \varepsilon +O(\varepsilon ^2)$, for some $c_3>0$
(the $e_i$'s are the standard unit vectors). 
Thus, we have to consider cross-polytopes $C_i$ contained in $K^*$, and
$C_o$ containing $K^*$, with centres at $o$. 
Of course, it suffices to show 
\begin{equation}
\label{intr}
|C_i|/|K^*| \le 1-c_4 \varepsilon +O(\varepsilon ^2), 
{\text{\,\,\,\,for some\,\,\,\,}}
c_4>0\,.
\end{equation}
We may assume that ${\text{vert}}\,C_i \subset {\text{vert}}\,K^*$.
Here ${\text{vert}}\,K^*$ consists of $\pm e_i$, for $1 \le i \le d$,
and still four vertices, close to $(\pm e_1 \pm e_2)/2$.
If for some $i \ge 3$ we have $\pm e_i \not\in C_i$, then $|C_i|=0$.
If $\pm e_1, \pm e_2 \in C_i$, then (\ref{intr}) holds. Otherwise, e.g., 
$\pm (1/2,1/2) \in {\text{vert}}\,C_i$, and either e.g.
$\pm e_1 \in {\text{vert}}
\,C_i$, or $\pm (1/2,-1/2)\in {\text{vert}}\,C_i$; in both cases $|C_i|/|K^*|=
1/2+ O( \varepsilon )$. So (\ref{intr}) is shown.

In a still more recent paper, namely Kim-Reisner \cite{KR}, there is 
proved the asymmetric variant of the theorem of
Nazarov-Petrov-Ryabogin-Zvavitch \cite{NPRZ}.
{\it{For $d \ge 2$ an integer there exist
$\varepsilon _d' >0$ and $c_d'>0$ with the following properties. 
If the Banach-Mazur distance of a convex body $K \subset {\mathbb{R}}^d$,
with}} $o \in {\text{int}}\,K$, {\it{from the class of
simplices is $1+ \varepsilon  \in (1, 1+ \varepsilon _d']$, 
then the volume product
$|K| \cdot |K^*|$ is at least $[(d+1)^{d+1}/(d!)^2](1+c_d' \varepsilon )$.
Again}}, also here {\it{the order of the error term is optimal}}. 
(An example is
obtained from a regular simplex of edge length $1$, and barycentre $o$, with 
small regular simplices of edge lengths $\varepsilon $ cut off at each
vertex. The argument showing optimality of the order of the error term
is like above.)

For general information about stability versions of geometric inequalities
cf. Groemer \cite{Gro}.

{\bf{Functional variants.}}
Also variants of the volume product problem have been treated.
E.g., {\it{functional forms
of the inverse Blaschke-Santal\'o inequality}} (i.e., of the lower estimate of
the volume product), cf. Meyer-Reisner \cite{MR98} 
(which states in p. 219 that a special case of its Theorem is 
the Mahler-Meyer theorem), {\it{functional forms of the Blaschke-Santal\'o
inequality}}, cf. 
Fradelizi-Meyer \cite{FM07} (which states in pp. 386-387, 393-394 that its
results imply the Blaschke-Santal\'o theorem --- with the equality case for
$o$-symmetry) and Artstein-Avidan-Klartag-Milman \cite{A-AKM} (which states
in p. 37 that its results imply the Blaschke-Santal\'o theorem, with the case
of equality --- however, this holds, strictly speaking, only in the
$o$-symmetric case, cf. this introduction, the second paragraph following this
paragraph). 

Functional variants are of different natures. 
E.g., the case of
{\it{``fractional dimension''}}, cf. Fradelizi-Meyer, \cite{FM}. Also,
convex bodies can be generalized to {\it{log-concave functions}}, 
i.e.,
functions ${\mathbb R}^d \to [0, \infty )$, whose logarithm is concave. 
To a convex body $K \subset {\mathbb R}^d$,
with $o \in {\text{int}}\,K$, one has to associate the log-concave function
$\exp (- \| x \| ^2_K /2)$, where $\| \cdot \| _K$ is the asymmetric norm 
with
unit ball $K$. 
Then $V(K)={\text{const}}_d 
\cdot \int _{{\mathbb R}^d} f(x)\,dx$, so 
here the integral on the 
right hand side is the proper substitute of $V(K)$.
Moreover, the polarity between $K$ and $K^*$ goes over to the following. If
we take the negative logarithms of two log-concave functions $f$ and $f^*$, 
then they are
the Legendre transforms of each other. The {\it{Legendre transform}} of a
function $\varphi : {\mathbb R}^d \to [ -\infty , \infty ]$ is 
${\cal L} \varphi : {\mathbb R}^d \to [ -\infty , \infty ]$, where  
$$
({\cal L} \varphi )(y)
:= \sup \{ \langle x,y \rangle -\varphi (x) \mid x \in {\mathbb R}^d \}\,.
$$ 
Thus, the subject of investigation is 
$$
\int _{{\mathbb R}^n} f(x)\,dx \cdot \int _{{\Bbb R}^n} f^*(x)\,dx\,,
$$ 
where one supposes
$$
\int _{{\mathbb R}^n} f(x)\,dx \in (0, \infty )\,.
$$
Cf. the nice expositions in Artstein-Avidan-Klartag-Milman \cite{A-AKM} and 
Klartag-Milman \cite{KM}. A straightforward calculation shows that this
product of integrals is invariant under non-singular
linear substitutions of the variable
$x$ (analogously as for the volume product of convex bodies), 
and under taking positive multiples of the function $f$.

Unfortunately, translations of convex bodies have no (good) generalizations
to log-concave functions. 
Thus, in place of a translation $K \mapsto K-x$, where $x \in {\text{int}}\,K$,
one considers an arbitrary translate of the function $f$ (i.e., $x \mapsto 
f(x-x_0)$). Then {\it{one proves the sharp upper bound $(2 \pi )^d$
for a suitable translate of the
original function}} $f$, {\it{cf. Artstein-Avidan, Klartag, Milman}} 
\cite{A-AKM}, 
{\it{Klartag-Milman}} \cite{KM}. {\it{Here, 
for even functions $f$, one may choose}} $x_0=o$ (as for
$o$-symmetric bodies one may consider $(K-o)^*$, see \S 1), 
and, more generally, {\it{in the
general case one can choose for $x_0$
the barycentre $\int xf(x)dx / \int f(x)dx$ of}} $f$ (as for general convex
bodies $K$ one may consider $[ K-b(K) ]^*$, 
see above in this introduction) cf. \cite{A-AKM}.
Of course, this does not concern the question of the
lower bound (as it is a minimum problem), 
but, in case of the upper bound, only the $o$-symmetric case of the
volume product problem generalizes this way to even log-concave functions.
(Observe that to translates of
convex bodies there do not correspond translates of the respective functions.
Already for $d=1$, with unit ball $[-1,1]$, to translates of $[-1,1]$,
by some $c
\in (-1,1)$, there correspond the functions ${\text{exp}}\, [  
-x^2/\left( 2(1+{\text{sg}}\,x \cdot c) ^2 \right) ]$, whose graphs have no
vertical axis of symmetry, while translates of the corresponding function 
have them. This shows that this functional variant of the Blaschke-Santal\'o
inequality generalizes the case of convex bodies
only in the $o$-symmetric case.)

{\it{For the upper bound, in the even case, the functional variant}}
(i.e., $\int _{{\mathbb R}^d} fdx \cdot \int _{{\mathbb R}^d} f^*dx \le (2 \pi
)^d$, cf., Ball \cite{Bal}, \cite{A-AKM}, Fradelizi-Meyer \cite{FM07})
immediately {\it{implies
the $o$-symmetric case of the volume product problem}}: namely, the extremal
even functions (up to constant factors)
are ones derived from $o$-symmetric convex bodies, more exactly, from 
$o$-symmetric
ellipsoids. 

{\it{For
the lower bound, the functional variant, i.e., $\int _{{\mathbb R}^d} fdx 
\cdot \int _{{\mathbb R}^d} f^*dx \ge \left( \pi / (2e) +o(1) \right) ^d$
for the case of even functions, and $\int _{{\mathbb R}^d} fdx 
\cdot \int _{{\mathbb R}^d} f^*dx \ge \left( \pi / (4e) +o(1) \right) ^d$
for the case of general functions, are 
proved in Fradelizi-Meyer}} \cite{FM08a}, {\it{Theorem 7.}}
(We remark that \cite{FM08a} Theorem 7 gave only the lower bound const$^d$,
using Bourgain-Milman's theorem, \cite{BoMi}. However, if we rather use the
currently best lower bounds for the $o$-symmetric, and the general cases, 
due to Kuperberg, G., \cite{K08}, both cited above, then, by the proof of
\cite{FM08a} Theorem 7,
we obtain the above given lower estimates.) 

The conjectured minima, for the even, or the general case, are  
$4^d$, or $e^d$, respectively, 
cf. Fradelizi-Meyer \cite{FM08b}. The minimizing functions
are conjectured to be, in a suitable system of coordinates, with the
origin at $o$, 
the following ones, cf. \cite{FM08b}, Conjectures $(1')$, $(2')$ (in the
general case this is stated there not so explicitly).
For the even case, 
$f(x_1,...,x_d)=$ const$\,\cdot \,e^{-\| (x_1,...,x_k) \| _{K_1}} 
\chi _{K_2} (x_{k+1},...,x_d)$, where $0 \le k
\le d$, and $K_1
\subset {\mathbb R}^k$ and $K_2 \subset {\mathbb R}^{d-k}$ are 
Hanner-Hansen-Lima
bodies, and $\chi $ denotes characteristic function. (The constant factor was
missing in \cite{FM08b}.)
For the general case, 
$f(x_1,...,x_d)=$ const$\,\cdot \prod _{i=1} ^d e^{-x_i} \chi _{[-1, \infty )}
(x_i)$. 
(We remark that \cite{FM08b} considered actually the class of all
translates of a function, and considered the infimum of the product
$\int _{{\mathbb R}^d} fdx \cdot \int _{{\mathbb R}^d} f^*dx$ on this class.
However, a straightforward calculation gives that the minimum of this product
on the class of all translates of the function
const$\,\cdot \prod _{i=1} ^d e^{-x_i} \chi _{[-1, \infty )}
(x_i)$
is attained, when we translate this function with the $0$ vector.)
Still we note that the quotients of the above lower
estimates and the conjectured
minima are $\left( \pi / (8e) +o(1) \right) ^d$, in the even case, and
$\left( \pi / (4e^2)+o(1) \right) ^d$, in the general case.

{\it{In the case of unconditional 
functions}} (i.e., $f(x_1,...,x_n)=f(|x_1|,...,|x_n|)$), 
{\it{the sharp lower bound, namely $4^d$, is known, cf. Fradelizi-Meyer}} 
\cite{FM08a}, \cite{FM08b}, {\it{with
the case of equality characterized in Fradelizi-Gordon-Meyer-Reisner}}
\cite{FGMR}: {\it{these are just the above given conjectured extremal
functions for the $o$-symmetric case.}}
Observe that {\it{this includes the case $d=1$, even functions}}.

{\it{For general functions, Fradelizi-Meyer}} \cite{FM08b} {\it{proved the
conjectured lower bound $e^d$ 
for functions $f$ being translates of functions $g$, that
vanish outside $[0, \infty )^d$, and are monotonically non-increasing in
each variable on $[0, \infty )^d$. 
They characterized the case of equality: these are just those of 
the above given
conjectured extremal functions for the general case, that verify the
hypotheses of their theorem.}}
Also {\it{for $d=1$, for general $f$, one has the sharp lower bound $e$,
cf. Fradelizi-Meyer}} \cite{FM08b}, {\it{who 
also determined the cases of equality: these are just the above given
conjectured extremal functions for the general case.}}

Still we note that the conjectures in ${\mathbb R}^d$ 
about the lower bound for the
functional variant, for the even, or the general case (i.e., $4^d$, or $e^d$), 
would imply the
conjectures about the lower bound for the volume product, in the
$o$-symmetric case, or in the case $o \in
$ int\,$K$, in ${\mathbb R}^d$, or ${\mathbb R}^{d-1}$
(that is, $4^d/d!$, or $d^d/[(d-1)!]^2]$), 
respectively, cf. \cite{FM08b}. 
However, {\it{the
conjecture about the lower bound for the
functional variant, for the even, or the general case,
for all $d$, is 
equivalent to the conjecture of Mahler-Guggenheimer-Saint Raymond, 
or of Mahler,
for all $d$, respectively, cf.}} \cite{FM08b}.

{\bf{Miscellaneous}} (added 5. April 2013)
\cite{Stancu} and 
\cite{RSchW} proved that the volume product $|\left( K-s(K) \right) ^*| \cdot
|K|$, or $|K^*| \cdot |K|$ can be (locally) 
minimal only if the generalized Gauss
curvature of $K$
is a.e. $0$, for the general, or $o$-symmetric case, respectively.
\cite{Kim} proved the following: if $K$ is $o$-symmetric, and $(1+\varepsilon
)$-close to some Hanner-Hansen-Lima body, in Banach-Mazur distance
(for suitable $\varepsilon >0$), then Mahler's conjecture holds for
$K$, and there is also a stability variant of this statement.
\cite{KimZva} proved a stability variant of the statement that among
unconditional bodies Mahler's conjecture is true. They also proved
the following: among $o$-symmetric convex bodies, $(1+
\varepsilon )$-close to unconditional convex bodies, in Banach-Mazur distance
(for suitable $\varepsilon >0$), the minimum of the volume
product is attained exactly for the Hanner-Hansen-Lima bodies, and they proved
also a certain
stability variant of this statement. \cite{BaBoFr} proved a stability
variant of functional forms of the Blaschke-Santal\'o inequality. \cite{GiPaVr}
also proved a version of the Bourgain-Milman theorem. A small survey is given
in \cite{Makai}.


\section{Main statements}

For stability versions of the Mahler-Reisner, Mahler-Meyer, and
Eggleston theorems, we prove the following theorems. 


As mentioned in the introduction, the following theorem was obtained also,
independently, by B\"or\"oczky, K. J.-Hug \cite{BH}, even for zonoids in
${\R}^d$ (the constants in \cite{BH} are unspecified, 
and stability of the centre of polarity is not investigated in \cite{BH}).

\begin{Thm}
\label{MahlerReisnerstab}
Let $K$ be a centrally symmetric convex body in $\R^2$ with $o \in 
{\text{\rm{int}}}\,K$ and $P$ a parallelogram, and
$$
|K|\cdot|K^*|\leq (1+\varepsilon)\cdot 8\,, 
{\text{\,\,\,\,with\,\,\,\,}} \varepsilon >0 \,.
$$
Then $\delta_{BM} (K,P) \leq 1+200\varepsilon $. Moreover, let $x \in 
{\mathbb R}^2$ and $\lambda _i >0$, and let $P$ be a parallelogram such that
$\lambda _1 P + x \subset K \subset \lambda _2 P + x$, and 
$\lambda _2 / \lambda _1 \le 1+200\varepsilon < 2 $. 
Then, in the Euclidean norm, for which
$[(\lambda _1 + \lambda _2 )/2]P$ is a square of diameter $1$,
we have that the distance of the
centre of 
$[(\lambda _1 + \lambda _2)/2] P + x$ 
from $o$ is at most $336 \cdot \sqrt{ \varepsilon }$. 
\end{Thm}


\begin{Thm}
\label{MahlerMeyerstab}
Let $K$ be a convex body in $\R^2$ with $o \in 
{\text{\rm{int}}}\,K$ and $T$ a triangle, and
$$
|K|\cdot|K^*|\leq (1+\varepsilon)\cdot 27/4\,,  
{\text{\,\,\,\,with\,\,\,\,}} \varepsilon >0 \,.
$$
Then 
$\delta_{BM}(K,T)\leq 1+900 \varepsilon $. 
Moreover, let $x \in 
{\mathbb R}^2$ and $\lambda _i >0$, and let $T$ be a triangle such that
$\lambda _1 T + x \subset K \subset \lambda _2 T + x $, and 
$\lambda _2 /
\lambda _1 \le 1+900 \varepsilon < 4$. Then, in the Euclidean norm, for which
$[(\lambda _1 + \lambda _2)/2]T$ is a regular triangle of side $1$, 
we have that the distance of the
centre of $[(\lambda _1 + \lambda _2)/2]T + x$ 
from $o$ is at most $917 \cdot \sqrt{ \varepsilon }$.
\end{Thm}


We note that, for $R_n$ a regular $n$-gon with centre $o$,
$$
|R_n|\cdot|R_n^*|=(n/2)\sin(2\pi/n)\cdot n\tan(\pi/n)=n^2\sin^2(\pi/n).
$$
We prove the following generalization of
the Mahler-Reisner and Mahler-Meyer theorems.


\begin{Thm}
\label{KiKo}
Let $K_i$ and $K_o$ be regular $n$-gons, $n\geq 3$,
and let each vertex of $K_i$ lie on a side of $K_o$, and
hence $K_i$ and $K_o$ have a common centroid $z$. If $K_i\subset K\subset 
K_o$ for
a planar convex body $K$ with $o\in{\rm int}\,K$, then
$$
|K|\cdot|K^*|\geq n^2\sin^2(\pi/n),
$$
with equality if and only if $o=z$, and either $K=K_i$, or $K=K_o$.
\end{Thm}


Let us show how Theorem~\ref{KiKo} yields the Mahler-Reisner and Mahler-Meyer 
theorems.
For the $o$-symmetric case, one considers
an ($o$-symmetric) parallelogram $P$ of maximal area contained in $K$.
Applying a linear map, we may assume that $P$ is a square.
Now the Mahler-Reisner theorem follows as $K\subset Q$ for the square $Q$ 
satisfying
that the midpoints of its sides are the vertices of $P$.

For the Mahler-Meyer theorem, 
let $T$ be a triangle of maximal area contained in $K$.
Applying a linear map, we may assume that $T$ is regular,
and let $S$ be the regular triangle satisfying
that the midpoints of the sides of $S$ are the vertices of $T$.
Since $K\subset S$,
Theorem~\ref{KiKo} yields the Mahler-Meyer theorem.


\vskip.3cm

Another consequence of Theorem~\ref{KiKo} is the following.


\begin{Cor}
\label{nfold}
If a convex body $K$ in $\R^2$ has $n$-fold rotational
symmetry about $o$, where $n\geq 3$, then
$$
|K|\cdot|K^*|\geq n^2\sin^2(\pi/n),
$$
with equality if and only if $K$ is a regular $n$-gon.
\end{Cor}


To prove Corollary~\ref{nfold} based on Theorem~\ref{KiKo}, one just
chooses a point $x\in\partial K$ that is the farthest from $o$,
and $K_i$ is the inscribed regular $n$-gon, of centre $o$, such that $x$ is
one of its vertices, and $K_o$ is the regular $n$-gon such that
the midpoints of the sides of $K_o$ are the vertices of $K_i$.


\begin{Thm}
\label{nfoldstab}
Let $n \ge 3$ be an integer, let $K$ be an $n$-fold 
rotationally symmetric convex body in $\R^2$ with $o \in 
{\text{\rm{int}}}\,K$ and $R_n$ a regular $n$-gon, 
and let
$$
|K|\cdot|K^*|\leq (1+\varepsilon)\cdot n^2\sin^2(\pi/n)\,, 
{\text{\,\,\,\,with\,\,\,\,}}
\varepsilon >0 \,.
$$
Then $\delta_{BM}^s (K,R_n)
\leq 1+18 \varepsilon $. Moreover, let
$x \in 
{\mathbb R}^2$ and $\lambda _i >0$, and let $R_n$ be 
a regular $n$-gon such that
$\lambda _1 R_n + x \subset K \subset \lambda _2 R_n + x$, and 
$\lambda _2 / \lambda _1 \le 1+18 \varepsilon < 1/\cos (\pi/ n) \le
2$. Then, in the Euclidean norm, for which
$[(\lambda _1 + \lambda _2)/2]R_n $ is a regular $n$-gon of diameter $1$, 
we have that the distance of the
centre of $[(\lambda _1 + \lambda _2)/2]R_n + x$ 
from $o$ is at most $263 \cdot \sqrt{ \varepsilon }$.
\end{Thm}


The following theorem proves the conjecture mentioned in \S 2, concerning
the exact error term in the stability variant of the Zhang projection body
inequality, for the planar case.


\begin{Thm}
\label{Egglestonstab}
Let $K$ be a convex body in $\R^2$ with
$$
|K| \cdot | \left( (K-K)/2 \right) ^*| \leq (1+\varepsilon)\cdot 6 \,,  
{\text{\,\,\,\,with\,\,\,\,}} \varepsilon >0 \,.
$$
Then $\delta_{BM}(K,T) \leq 1+87 \varepsilon $.
\end{Thm}


{\bf{Remark.}}
In Theorems 1, 2 we restricted ourselves to the case $\lambda _2 / \lambda _1
<2$, or, $\lambda _2 / \lambda _1 <4$, respectively. This we did since 
$\lambda _2 / \lambda _1 =2$, or $\lambda _2 / \lambda _1 =4$ is no restriction
at all for the body $K$. Namely, if $K_i \subset K$ is an $o$-symmetric
inscribed parallelogram of maximal area (for $K$ $o$-symmetric), or an
inscribed triangle of maximal area (for $o \in $ int\,$K$), 
then $K \subset K_o$, where $K_o$ is a
parallelogram, or triangle, with side midpoints at the vertices of $K_i$.
If $K'$ is a parallelogram, or triangle, with side midpoints at the vertices
of $K_o$, then $K_i$ and $K'$ are homothetic, with homothety ratio $2$, or $4$,
respectively, and $K_i \subset K \subset K_o \subset K'$. (For Theorem 5 the
analogous inequality would be $\lambda _2 / \lambda _1 \le 1/ \cos ^2 (\pi /n)
<4$, but this does not suffice to prove (\ref{epsilon}), 
in the proof of Theorem \ref{nfoldstab}, with positive right
hand side.)


{\bf Example.} {\bf{1.}}
We show that the stability statements in  
Theorems~\ref{MahlerReisnerstab}, ~\ref{MahlerMeyerstab}, 
~\ref{nfoldstab}, ~\ref{Egglestonstab}, concerning the bodies, 
are of the exact order of magnitude. For this, let
the regular $n$-gon $R_n$ be inscribed in the unit circle $U$ about $o$, 
and let
us define $K_n$ as the convex polygon with vertices the vertices of $R_n$, and
$1+ \varepsilon $ times the side-midpoints of $R_n$, where $\varepsilon \in
(0, 1/ \cos ( \pi /n)]$ (thus $K_n \subset U$). 
Then $|K_n| \cdot |(K_n)^*|=n^2 \sin ^2 ( \pi /n) + n^2 \sin ^2 ( \pi /n)
\cdot \left( \varepsilon -\varepsilon ^2 \cot ^2 ( \pi /n) \right)/(1+
\varepsilon )$. Letting $n=3$, we have $|K_3| \cdot | \left( (K_3-K_3)/2
\right) ^*| = 6 \cdot (9+15 \varepsilon + 3 \varepsilon ^2 - 3 \varepsilon ^3)
/(3+ \varepsilon )^2$.
Clearly, $\delta _{BM}^s(R_n,K_n) \le 1+ \varepsilon
$. On the other hand, for suitable $A$ and $x$, we have 
$ \lambda _1 AR_n \subset K_n \subset \lambda _2 AR_n +x $ and
$\delta _{BM}(R_n,K_n)^2=( \lambda _2 / \lambda _1 )^2 \ge |K_n|/|\lambda _1
AR_n| \ge (1 + \varepsilon ) |R_n|/|R_n|$ (at the last step we have used that 
$ \lambda _1 AR_n \subset U$ is a convex $n$-gon, hence $|\lambda _1 AR_n| \le
|R_n|$, similarly as in the end of the introduction, at the proof of the
optimality of the order of the error term). 
Hence, $\delta _{BM}(R_n,K_n) \ge {\sqrt{1+ \varepsilon }} $. (For
Theorems~\ref{MahlerReisnerstab}, ~\ref{MahlerMeyerstab} we use the cases
$n=4,\,\,3$.) 

{\bf{2.}}  
For the stability of the centre of polarity (for Theorems 
\ref{MahlerReisnerstab}, \ref{MahlerMeyerstab}, \ref{nfoldstab}), we proceed
analogously to \cite{KR}, Proposition 2.
An example is a
regular $n$-gon $K$ of centre $o$, and diameter $1$ (with $\lambda _i=1$). 
We use the well-known formula (\ref{derivatives}) from the proof of 
Lemma~\ref{centrestab}, which comes in the sequel, for $d=2$. 
The inradius of $K$ is at least $1/(2
\sqrt{3})$. We let $\| x \| \le 1/(4 \sqrt{3})$, and estimate $(\partial /
\partial x_2)^2 | (K-x)^* |$ from above
by replacing, in the inequality in (\ref{derivatives}), $h_K(u)$ by
$1/(2 \sqrt{3})$, and then $(1- \langle u,x \rangle )^{-4}$ by 
$ \left( 1/(4 \sqrt{3}) \right) ^{-4}$. Then, using still
$\int _{S^1}u_2^2du=\pi $, we get 
$$
(\partial / \partial x_2) | (K-x)^* | =0 {\text{\,\,\,\,and\,\,\,\,}} 
(\partial / \partial x_2)^2 | (K-x)^* | \le 2^8 \cdot 3^3 \cdot \pi \,.
$$
By ${\text{diam}}\,K =1$ we have $|K| \le \pi/4$. Thus we get
$(\partial / \partial x_2)^2 (|K| \cdot | (K-x)^* |) \le 2^6 \cdot 3^3 \cdot
\pi ^2$, and the
analogues of these formulas hold for the first and 
second directional derivatives in any
direction. Thus, for $|K| \cdot |(K-x|^*| \ge (1+\varepsilon ) \cdot n^2 \sin
^2 ( \pi /n )$, we have
$$
\varepsilon \cdot 27/4 \le \varepsilon \cdot n^2 \sin ^2 ( \pi /n ) \le
|K| \cdot | (K-x)^* | - |K| \cdot |K^*| \le 2^5 \cdot 3^3 \cdot \pi ^2 \| x \|
^2\,,
$$
hence, for any $x$ --- i.e., without the restriction $\| x \| \le 1/(4
\sqrt{3})$ --- we have
 $$ 
\| x \| \ge \sqrt{ \varepsilon } \cdot \sqrt{2} 
/ (16 \pi ) {\text{\,\,\,\,or\,\,\,\,}} 
\| x \| \ge 1/(4 \sqrt{3})\,.
$$
Then the first one of these inequalities holds, 
if $\varepsilon \in (0, \varepsilon _0]$,
where $\sqrt{ \varepsilon _0} \cdot \sqrt{2} /(16 \pi ) =1/(4 \sqrt{3})$,
i.e., for $\varepsilon _0= 8 \pi ^2 /3$. $ \blacksquare $

\hskip0.0cm


In a forthcoming paper, by the first two named authors of
this paper, 
we will show that, for convex $n$-gons $K$, 
the product $|K| \cdot |  [ K-s(K) ] ^* | $ is
maximal exactly for the affine regular $n$-gons. We remark that the
$o$-symmetric case of this statement is obtained, independently, also by
Tabachnikov, in a more general form, namely for star-polygons, \cite{Ta},
Theorem 2.
Moreover, the general case is obtained independently, also by the last two
named authors of this paper; see \cite{MR}.
Further, in the above mentioned forthcoming paper, we will give stability
estimates for the Blaschke-Santal\'o inequality in the plane, 
for the $o$-symmetric case.
Here the deviation from the ellipses will be
measured by the quotient of the areas of the convex body, and the maximal area
inscribed/minimal area circumscribed ellipse of the convex body, and the order
of the error term will be optimal.
If any of these ellipses is the unit circle about $o$, then even the
arithmetic mean of the areas of the body and the polar body is at most $ \pi $.


\section{Proof of Theorem~\ref{KiKo}}
\label{secKiKo}

First we prove a lower bound for the volume product
in sectors. The idea of giving lower bounds in sectors separately, and then
using the arithmetic-geometric mean inequality, is due to Saint Raymond
\cite{SR}, proof
of Th\'eor\`eme 28. There it is also noted that this approach settles the
two-dimensional $o$-symmetric case. Our proofs of our Theorems 
~\ref{MahlerReisnerstab}, ~\ref{MahlerMeyerstab}, ~\ref{KiKo}, 
~\ref{nfoldstab} all use this idea.


\vskip.3cm

The particular case
$u=u^*=(0,1)$, and $v=v^*=(1,0)$ of our following lemma reduces to 
the two-dimensional case of \cite{SR}, Th\'eor\`eme 28.


\begin{Lem}
\label{section}
Let $K$ be a planar convex body with
$o\in{\rm int}\,K$. Let, for some linearly independent $u,v\in \partial K$,
and linearly independent $u^*,v^*\in\partial K^*$,
the supporting lines to $K$ with exterior normals $u^*$ and $v^*$
intersect $K$, e.g., at $u$ and $v$, respectively, and intersect
each other at $p \in \R^2$, where $[p,o] \cap [u,v] \neq \emptyset $.
Furthermore, let
the supporting lines to $K^*$ with exterior normals $u$ and $v$
intersect $K^*$, e.g., at $u^*$ and $v^*$, respectively, and intersect
each other at
$p^* \in \R^2$ with $[p^*,o] \cap [u^*,v^*] \neq \emptyset $.
Then, for $C=K\cap [o,u,v,p]$ and
$C^*=K^*\cap [o,u^*,v^*,p^*]$, we have
$$
|C|\cdot |C^*|
\geq |[o,u,v,p]|\cdot |[o,u^*,v^*]|=|[o,u,v]|\cdot |[o,u^*,v^*,p^*]|
\,,
$$
with equality if and only if either
$C=[o,u,v]$ or
$C=[o,u,v,p]$.
\end{Lem}


{\bf Remark.} 
We may assume $C \ne [o,u,v]$. Then, for
$p= \lambda u + \mu v$ and $\lambda , \mu >0$,
we have $\lambda+\mu>1$
and $p^*=\mu u^*+\lambda v^*$. We choose a coordinate system,
assuming 
$$
u=(1,0),\,\, {\rm{and}}\,\,v=(0,1)\,.
$$
Then
$$
p=(\lambda,\mu), \,\, p^*=(1,1),\,\, 
u^*=\left( 1,(1-\lambda )/\mu \right) ,\,\, 
v^*=\left( (1-\mu )/\lambda ,1 \right) \,,
$$
and
$$
\begin{cases}
|[o,u,v,p]|\cdot |[o,u^*,v^*]|= |[o,u,v]|\cdot |[o,u^*,v^*,p^*]| = \\
(\lambda +\mu )(\lambda + \mu -1)/(4 \lambda \mu )
=(2-\langle u,v^*\rangle-\langle u^*,v\rangle)/4\,.
\end{cases}
$$

\vskip.2cm


First we show that Mahler's original proofs (\cite{Mah38}) 
yield our lemma, except the case of equality. 

\vskip.3cm


{\bf First proof.}
We exclude $C=[o,u,v],[o,u,v,p]$.
Let $k \ge 0$ be an integer, and let us suppose that both $C$ and $C^*$ are
polygons, such that the total number 
of their vertices in ${\text{int}}\,[u,v,p]$, or
${\text{int}}\,[u^*,v^*,p^*]$, respectively, is at most $k$. 
(This case suffices to prove the
inequality.) Let $C, C^*$ realize the minimum under these hypotheses. If e.g.
$C$ has a vertex $c \in {\text{int}}\,[u,v,p]$, then we can move $c$ a bit,
parallel
to the diagonal connecting its neighbours, hence keeping $|C|$ fixed. 
Then, for $C^*$, the polar side line
will rotate about some of its points. Since the lines of the neighbours of 
this side intersect outside this side line, by some small rotation $|C^*|$
strictly decreases, a contradiction. Hence we have a situation as for $k=0$.

If $k=0$, then $C$ has a vertex $c$, e.g. in ${\text{relint}}\,[u,p]$, and then
$C=[o,u,v,c]$, since else $C^*$ would have a vertex in
${\text{int}}\,[u^*,v^*,p^*]$. Then $c=( \alpha
\lambda +1 - \alpha , \alpha \mu)$, where $ \alpha \in  (0,1)$, and
$|C| \cdot |C^*|=(1/4) \left( 1+( \lambda + \mu -1) \alpha \right) \cdot
[1-(1- \lambda /\mu )-\left( (1- \mu )/\lambda -1 \right) (1-\alpha \lambda )/
(1- \alpha + \alpha \mu)]$.
The fact that this is at least $(\lambda + \mu )(\lambda + \mu -1)/(4 \lambda
\mu )$ can be written, after multiplying with the product of the
denominators (each of them being positive), and rearranging (using the program
package GAP
\cite{GAP}), as
$\lambda (\lambda + \mu -1)^2 \cdot \alpha (1-\alpha ) \ge 0$.
$ \blacksquare $ 

\vskip.3cm


The second proof follows the lines of Meyer \cite{Me86}, proof of
Th\'eor\`eme I. 2 (more exactly, its two-dimensional case, that gives our
lemma for $u=u^*=(0,1)$, and $v=v^*=(1,0)$).

\vskip.3cm

{\bf Second proof.} We have
$$
1=\langle u^*,u\rangle=\langle u^*,p\rangle
=\langle v^*,p\rangle=\langle v^*,v\rangle=\langle u,p^*\rangle
=\langle v,p^*\rangle.
$$
For $x\in K\cap[p,u,v]$, the sum of the heights of the triangles $[o,u,v]$ and 
$[x,u,v]$, belonging to their common side $[u,v]$,
is $\langle p^*,x\rangle/\|p^*\|$. Thus
the vectors $w:=[\|u^*-v^*\| /(2\|p\| )]p$ and 
$w^*:=[\|u-v\| /(2\|p^*\| )]p^*$
satisfy
\begin{eqnarray}
\label{Kx}
 |C|  &\geq& |[o,u,v,x]|=\langle w^*,x\rangle
 \mbox{ \ for \ $x\in K\cap[u,v,p]$}, \mbox{ \ and} \\
 \label{K*y}
 |C^*| &\geq&
 |[o,u^*,v^*,x^*]|=\langle w,x^*\rangle
 \mbox{ \ for \ $x^*\in K^*\cap[u^*,v^*,p^*]$}.
\end{eqnarray}
Since $\langle w^*,p\rangle=|[o,u,v,p]|\geq |C|$,
and $\langle w^*,x\rangle<\langle w^*,u\rangle$
for $x\in K\backslash[p,u,v]$,
 we have $\widetilde{w}^*:=|C|^{-1}w^*\in K^*\cap[u^*,v^*,p^*]$,
and analogously $\widetilde{w}:=|C^*|^{-1}w\in K\cap[u,v,p]$.
It follows by applying (\ref{Kx}) to $x=\widetilde{w}$, that
$$
\begin{cases}
|C|\cdot |C^*|\geq \langle w^*,|C^*|\widetilde{w}\rangle=
\langle w^*,w\rangle= \\
\langle w^*,p\rangle\cdot\mbox{$\|u^*-v^*\|/(2\|p\|)$}=
|[o,u,v,p]|\cdot |[o,u^*,v^*]|.
\end{cases}
$$
We also have $\langle w^*,w\rangle=|[o,u,v]|\cdot |[o,u^*,v^*,p^*]|$
by the remark following the statement of this Lemma,
hence we have equality in the Lemma if $C=[o,u,v]$ or
$C^*=[o,u,v,p]$.

Assume that equality holds in Lemma~\ref{section}. It follows by
(\ref{Kx}) and (\ref{K*y}) that
$$
C=[o,u,v,\widetilde{w}] \mbox{ \ and \ }
C^*=[o,u^*,v^*,\widetilde{w}^*].
$$
In particular $C^*$ has vertices $a^*$ and $b^*$ satisfying
$$
\langle a^*,u\rangle=\langle a^*,\widetilde{w}\rangle=1
\mbox{ \ and \ }
\langle b^*,v\rangle=\langle b^*,\widetilde{w}\rangle=1.
$$
Checking the vertices of $C^*$, we have only two choices.
Either $a^*=u^*$ and $b^*=v^*$, and hence $C=[o,u,v,p]$,
or $a^*=b^*=\widetilde{w}^*$, and hence $C=[o,u,v]$.
$\blacksquare $


\vskip.3cm

The third proof will use an idea of Behrend, \cite{Be}, proof of (77),
pp. 739-740, and of (112), pp. 746-747.
Its idea, intuitively, is the following. 
``If $C$ is close to $[o,u,v]$, then $C^*$ is close to  
$[o,u^*,v^*,p^*]$, hence $|C^*|$ will be a lot greater than $|[o,u^*,v^*]|$. 
On the other hand, if $C$ is close to $[o,u,v,p]$, then $|C|$ 
will be a lot greater than $|[o,u,v]|$.'' 

\vskip.3cm

{\bf Third proof.}
Using the notations of the second proof, we have 
$$
|C| \ge |[o,u,v,x]|,
$$
where now $x$ is a point of $C \cap [u,v,p]$, that is 
farthest from $(p^*)^{-1}(1)$, which line passes
through $u,v$. Then there is a supporting line $(x^*)^{-1}(1)$ at $x$ to $K$,
parallel to $(p^*)^{-1}(1)$. Then
$$
|C^*| \ge |[o,u^*,v^*,x^*]|\,,
$$
so,
$$
|C| \cdot |C^*| \ge |[o,u,v,x]| \cdot |[o,u^*,v^*,x^*]|\,.
$$
Observe that, if $x$ varies in $[u,v,p]$, then $|[o,u,v,x]|$ 
is proportional to 
${\text{dist}} \left( o, \right.$ $\left. (x^*)^{-1}(1) \right) 
=1/\|x^* \|$. Simultaneously, $x^*$ varies in
$[o,p^*] \cap [u^*,v^*,p^*]$, hence $|[o,u^*,v^*,x^*]|$ is proportional to 
$\| x^* \|
$. Hence, $|[o,u,v,x]| \cdot |[o,u^*,v^*,x^*]|$ does not depend on $x$, so
has the same value, as for $x \in [u,v]$, and for $x=p$. 

For the case of equality we have $C=[o,u,v,x]$ and
$C^*=[o,u^*,v^*,x^*]$. We exclude $x \in [u,v]$ and $x=p$. Then $x^*$ varies
in ${\text{relint}}\,(C^* \cap [o,p^*])$, and we get a contradiction as in the
second proof.
$\blacksquare $

\vskip.3cm


{\bf Proof of Theorem~\ref{KiKo}.}
We may assume that $o$ is the Santal\'o point of $K$.
First we show that $o\in {\rm int}\,K_i$.

We note that as
the origin is the centroid of $K^*$, there exists
no line $l$, with $o\in l$, and bounding the
half planes $l^-$ and $l^+$, such that
the reflected image of $K\cap l^-$ through the line $l$
is strictly contained in $K\cap l^+$.
If $n\geq 4$ then the angles of a regular $n$-gon are at least $\pi/2$,
thus $o\in {\rm int}\,K_i$ by the property of the Santal\'o point above.

If $n=3$ then we may assume that $K$ is not a parallelogram.
In this case for each triangle $S$ cut off by a side $s$
of $K_i$ from $K_o$, there is a linear transformation $A$
such that the reflected image of $A(S)$
through the line $A(s)$ is strictly contained in $A(K_i)$ (here we use that $K$
is not a parallelogram). Therefore
the property of the Santal\'o point above, applied to $A(K)$,
yields $o\in {\rm int}\,K_i$.

When indexing the vertices of an $n$-gon, we identify
vertices with indices $j$ and $j\pm n$.
Let $x_1,\ldots,x_n$, and $y_1,\ldots,y_n$ denote
the vertices of $K_i$ and $K_o$ in counterclockwise order, and
$x^*_1,\ldots,x^*_n$, and $y^*_1,\ldots,y^*_n$ denote
the vertices of $K_i^*$ and $K_o^*$, respectively, so
that, for $j=1,\ldots,n$, we have $x_j\in[y_j,y_{j+1}]$, and
$$
1=\langle x_j^*,x_{j-1}\rangle=\langle x_j^*,x_j\rangle
=\langle y_j^*,y_{j+1}\rangle=\langle y_j^*,y_j\rangle.
$$
In particular, $y_j^*\in[x_j^*,x_{j+1}^*]$.
For $j=1,\ldots,n$,
let $C_j=K\cap[o,x_{j-1},x_j,y_j]$ and
 $C_j^*=K^*\cap[o,y_{j-1}^*,y_j^*,x_j^*]$.
 Therefore Lemma~\ref{section} yields that
 \begin{equation}
 \label{CjC*j}
|C_j|\cdot |C_j^*|\geq |[o,x_{j-1},x_j,y_j]|\cdot
|[o,y_{j-1}^*,y_j^*]|,
 \end{equation}
 with equality if and only if $C_j=[o,x_{j-1},x_j,y_j]$
 or $C_j=[o,x_{j-1},x_j]$.

 By the $n$-fold rotational symmetry of $K_i$ and $K_o$
 about their common centre, there exist common distances
 $a=\|x_{j-1}-y_j\|$ and $b=\|x_j-y_j\|$ for $j=1,\ldots,n$,
 and hence $a+b$ is the side length of $K_o$.
 Since the distance
 of $o$ from the line $y_jy_{j+1}$ is
 $d_j:=\|y_j^*\|^{-1}$, for $j=1,\ldots,n$, it follows that
$$
|C_j|\cdot |C_j^*|=\frac{(ad_{j-1}+bd_j)\sin(2\pi/n)}{4d_{j-1}d_j}\,.
$$
Additionally, we have
$$
\frac{n(a+b)^2}{4\tan(\pi/n)}=|K_o|=\frac{(a+b)(d_1+\ldots+d_n)}2\,.
$$
We deduce by repeated applications of the inequality between
the (weighted) arithmetic and geometric means, that
\begin{eqnarray}
\label{sumCj}
|K|\cdot|K^*| &=&\left(\sum_{j=1}^n|C_j|\right) \cdot
\left(\sum_{j=1}^n|C_j^*|\right)
\geq n^2\left(\prod_{j=1}^n ( |C_j|\cdot |C_j^*| ) \right)^{1/n}\\
\nonumber
 &=& \frac{n^2\sin(2\pi/n)}4
 \left(\prod_{j=1}^n\frac{ad_{j-1}+bd_j}{d_{j-1}d_j}\right)^{1/n}  \\
 \label{sumproddj}
&\geq&
 \frac{n^2(a+b)\sin(2\pi/n)}{4}
 \left(\prod_{j=1}^n d_j\right)^{-1/n}\\
 \label{proddj}
 &\geq &\frac{n^3(a+b)\sin(2\pi/n)}{4\sum_{j=1}^n d_j}
=
 \frac{n^2\sin(2\pi/n)\tan(\pi/n)}{2}\,.
\end{eqnarray}

 Assume that equality holds in Theorem~\ref{KiKo}. It follows by
 (\ref{sumproddj}) and 
 (\ref{proddj}) that all $d_j$ are equal, thus $o$ is the
 common centre of $K_i$ and $K_o$. Further,
 all $C_j$ have the same area by (\ref{sumCj}).
 Therefore the equality conditions
 in (\ref{CjC*j}) imply that either $K=K_i$
 or $K=K_o$. 
$\blacksquare $

\vskip.3cm


{\bf{Remark.}}
In the particular case of Lemma~\ref{section}, when $C$ is an $n$-th part
of a convex body $K$ with $n$-fold rotational symmetry about $o$,   
we could have referred in the first proof to 
\cite{MR06}, 
to the so called ``shadow movement'' (although this is more involved than
the elementary proof of Mahler used above). 
That is, we have an $ln$-gon
$K=x_1...x_{ln}$ (where $l \ge 2$), 
having $n$-fold rotational symmetry about $o$.
The movement of the vertices $x_2, x_{2+l},...x_{2+(n-1)l}$, parallel to the
diagonals $x_1x_3$, etc., preserving the rotational symmetry, and giving a
polygon $K'$, of course does not determine a
shadow movement. However, we can move only $x_2$, in the above way, and 
this determines a shadow movement, giving a polygon $K''$. (More exactly: only
the points of $[x_1,x_2,x_3]$ are moved, in the direction of $x_1x_3$. At this
motion, 
the points of any chord, parallel to $x_1x_3$, 
are moved with the same velocity,
so that at any moment the moved chords constitute a triangle with vertices
$x_1,x_3$, and the translate of $x_2$). 
Then $|K|=|K'|=|K''|$, and 
$| (K')^* | =|K^*|+n (| (K'')^* | - |K^*|)$, so 
$| (K')^* | $ is a linear function of $| (K'')^* |$.
Moreover, $K''$ and $K$ are not affinely equivalent (consider the barycentres
of the subpolygons with vertices each $l$'th vertex of $K'',K$).


\section{Proofs of the stability theorems}

The main result in this section
is the following stability version
of Lemma~\ref{section}.


\begin{Lem}
\label{sectionstab}
Let $C,C^*,u,u^*,v,v^*,p,p^*$ be as in Lemma~\ref{section},
and let $p=\lambda u+\mu v$ for $\lambda ,\mu>0$.
If
$$
|C| \cdot |C^*|
\leq (1+\varepsilon) |[o,u,v,p]|\cdot |[o,u^*,v^*]|,
$$
for positive $\varepsilon < \min \, \{ \lambda ,\mu \}/(\lambda +\mu )$,
then for 
$\gamma :=3[(\lambda+\mu)/(\min \, \{ \lambda,\mu \})]
(1+\sqrt{\lambda +\mu })$,
$$
\begin{cases}
{\text{either\,\,\,\,}} C\subset(1+\gamma\varepsilon)[o,u,v]\,, \\  
{\text{or\,\,\,\,}} (1+\gamma\varepsilon)^{-1}[o,u,v,p]\subset C,
{\text{\,\,\,\,that is\,\,\,\,}} C^* \subset(1+\gamma\varepsilon)[o,u^*,v^*]\,.
\end{cases}
$$
\end{Lem}


{\bf First proof.} 
We may assume $C \ne [o,u,v]$.
We use the notations from the Remark after Lemma~\ref{section}, and from
the second proof of Lemma~\ref{section}.
We have $\widetilde{w}=tp$ and $\widetilde{w}^*=sp^*$ for some $t,s\in(0,1]$.
Since $\langle \widetilde{w},\widetilde{w}^*\rangle\leq 1$, we have
\begin{equation}
\label{ts}
ts(\lambda+\mu)\leq 1.
\end{equation}
Further, for $\tilde{u}^*:=\left(1,(1-t\lambda )/(t\mu )\right)$ and 
$\tilde{v}^*:=\left( (1-t\mu)/(t\lambda ),1\right)$, we have
$$
1=\langle \tilde{u}^*,u\rangle=\langle \tilde{u}^*,\widetilde{w}\rangle
=\langle \tilde{v}^*,v\rangle=\langle \tilde{v}^*,\widetilde{w}\rangle.
$$
It follows by the second proof of Lemma~\ref{section}, using the notations
$\widetilde{w},{\widetilde{w}}^*$ introduced there, that
\begin{eqnarray}
\label{uvinC}
[o,u,v,\widetilde{w}]\subset C 
\mbox{ \ and \ }
|C|\leq(1+\varepsilon)|[o,u,v,\widetilde{w}]|, \mbox{ \ and} \\
\label{u*v*inC*}
[o,u^*,v^*,\widetilde{w}^*]\subset C^* 
\mbox{ \ and \ }
|C^*|\leq(1+\varepsilon)|[o,u^*,v^*,\widetilde{w}^*]|.
\end{eqnarray}
It follows that if 
$\langle\tilde{u}^*,x\rangle\geq\langle\tilde{u}^*,u\rangle=1$
for $x\in C$ then
$$
|[x,u,\widetilde{w}]|\leq \varepsilon\cdot|[o,u,v,\widetilde{w}]|
=\varepsilon\cdot
\mbox{$ [ (\lambda+\mu )/\mu ] $}\cdot|[o,u,\widetilde{w}]|,
$$
and hence $\langle\tilde{u}^*,x\rangle\leq 
1+\varepsilon\cdot (\lambda+\mu)/\mu $.
For $\widetilde{\gamma}:=(\lambda+\mu)/\min \, \{\lambda,\mu\} $,
we deduce that 
$C\subset(1+\widetilde{\gamma}\cdot\varepsilon)[o,u,v,\widetilde{w}]$,
and hence $[o,u^*,\tilde{u}^*,v^*,\tilde{v}^*,\widetilde{w}^*]
\subset(1+\widetilde{\gamma}\cdot\varepsilon)C^*$ by polarity,
and analogously 
$C^*\subset(1+\widetilde{\gamma}\cdot\varepsilon)[o,u^*,v^*,\widetilde{w}^*]$.
Since $\varepsilon<\widetilde{\gamma}^{-1}$, we deduce
$$
[o,u^*,\tilde{u}^*,v^*,\tilde{v}^*,\widetilde{w}^*]
\subset(1+\widetilde{\gamma}\cdot\varepsilon)^2[o,u^*,v^*,\widetilde{w}^*]
\subset(1+3\widetilde{\gamma}\cdot\varepsilon)[o,u^*,v^*,\widetilde{w}^*].
$$
For $a:=(\lambda-s\lambda,s\lambda+\mu-1)$, we have 
$\langle a,v^*\rangle=\langle a,\widetilde{w}^*\rangle=s(\lambda+\mu-1)$, thus
\begin{eqnarray*}
1+3\widetilde{\gamma}\cdot\varepsilon&\geq&
\frac{\langle a,\tilde{v}^*\rangle}{\langle a,v^*\rangle}=
\frac{ts(\lambda+\mu-1)+(1-s)(1-t)}{ts(\lambda+\mu-1)}\\
&=&
1+\left(\frac1s-1\right)\left(\frac1t-1\right)\frac1{\lambda+\mu-1}.
\end{eqnarray*}
It follows by (\ref{ts}) that 
$$
{\text{either\,\,\,\,}} 1/s\geq\sqrt{\lambda+\mu}, 
{\text{\,\,\,\,or\,\,\,\,}} 
1/t\geq\sqrt{\lambda+\mu}\,.
$$ 
In the first
case, $3\widetilde{\gamma}\cdot(\lambda+\mu-1)/(\sqrt{\lambda+\mu}-1)=\gamma$ 
yields
$1/t \le 1+\gamma\varepsilon$, and
hence $(1+\gamma\varepsilon)^{-1}[o,u,v,p]\subset C$.
On the other hand, if $1/t\geq\sqrt{\lambda+\mu}$, then
a similar argument leads to
$(1+\gamma\varepsilon)^{-1}[o,u^*,v^*,p^*]\subset C^*$,
and hence $C\subset (1+\gamma\varepsilon)[o,u,v]$.
$\blacksquare $


\vskip.3cm

The second proof of Lemma~\ref{sectionstab}, where however
the constant $\gamma $ will be different, and
which iterates the construction in the proof of 
Behrend (\cite{Be}, proof of (77), pp. 739-740, and of (112), pp. 746-747)
will be broken up into two parts.


\begin{Lem}
\label{sectionstab2}
Under the hypotheses of Lemma~\ref{section}, and with
$p=\lambda u+\mu v$, for $\lambda,\mu>0$,
we have 
$$
|C| \cdot |C^*| \ge f( \lambda , \mu )
+g(\lambda ,\mu ) \alpha (1-\alpha )\,,
$$
where 
$$
f( \lambda , \mu ):=
(\lambda + \mu )(\lambda + \mu -1)/(4 \lambda \mu )\,,
$$ 
$$g( \lambda , \mu
):=
(1/4) \cdot  ( \lambda + \mu -1)^2 \cdot  \min \, \{ 1/[ \mu (1+\lambda /4
+\mu )], 1/ [ \lambda (1+\lambda + \mu /4 )], 1/(\lambda \mu) \} \,,
$$ 
$$
\alpha := \max \, \{ |[u,v,x]|/|[u,v,p]| \mid x \in C \cap [u,v,p] \} \in 
[0,1]\,.
$$
\end{Lem}


{\bf Proof.}
Again we use the notations from the Remark after Lemma \ref{section}.

We may suppose $\alpha \in (0,1)$.
Let $x =(x_1,x_2)\in C \setminus [o,u,v]$ 
realize $\alpha =\max |[u,v,x]|/|[u,v,p]$. We
write $C_i:=[o,u,v,x]$, and 
$C_o:=\{ ( \xi , \eta ) \in [o,u,v,p] \mid \xi + \eta 
\le x_1+x_2 \} $. Then $C_i \subset C \subset C_o$.
Let $x$ divide the chord of $[o,u,v,p]$, parallel to the line $uv$, and 
containing $x$ 
in the ratio $\beta : (1 - \beta )$, where $\beta \in [0,1]$, and
where the part of the chord with ratio
$\beta $ has an endpoint in $[u,p]$.

We iterate this construction. Let $y,z \in C$, and $\overline{y}, \overline{z} 
\in C_o$ lie on the other sides of the
lines $ux,vx$ than $o$, and let them realize $\max |[u,x,y]|$, or 
$\max |[v,x,z]|$ and
$\max |[u,x,\overline{y}]|$, or $\max |[v,x,\overline{z}]|$
under these conditions, respectively. 
We define $\gamma :=|[u,x,y]/|[u,x,\overline{y}]| \in
[0,1]$ and 
$\delta :=|[v,x,z]|/$
\newline
$|[v,x,\overline{z}]| \in [0,1]$. 
Let $C_i':=C_i \cup [u,x,y] \cup [v,x,z]$, and let
$C_o'$ be the intersection of $C_o$ and the support half-planes of $C$ at
$y,z$, with boundaries parallel to the lines $ux,vx$. Then $C_i \subset C_i'
\subset C \subset C_o' \subset C_o$. So for their polars (in the angular
domain $u^*ov^*$) we have $(C_o)^* \subset (C_o')^* \subset C^*$. Hence,
$$
\begin{cases}
|C| \cdot |C^*| \ge |C_i'| \cdot |(C_o')^*| \ge \\ 
|C_i| \cdot |(C_o)^*| +
|C_i' \setminus C_i| \cdot |(C_o)^*| + |C_i| \cdot |(C_o')^* \setminus
(C_o)^*| = \\
|C_i| \cdot |(C_o)^*| + (|T_y|+|T_z|) \cdot |(C_o)^*| +|C_i| \cdot
(|(T^*)_y|+|(T^*)_z|)\,.
\end{cases}
$$
Here $T_y:=[u,x,y]$ and 
$T_z:=[v,x,z]$, and the triangles $(T^*)_y$ and $(T^*)_z$ have
as their vertices the polars of the three first, or three last 
consecutive side lines of $C_o'$ in the open angular domain $u^*ov^*$,
taken in the positive orientation, respectively.

First we estimate $|T_y| \cdot |(C_o)^*| +|C_i| \cdot |(T^*)_y|$ from below.
We have 
$$
|C_i|=[1+(\lambda + \mu -1) \alpha ]/2\,, 
$$
$$
|(C_o)^*|=(1/2) \cdot [1/\left( 1+ (\lambda + \mu -1) \alpha \right) ] \cdot
( \lambda + \mu )( \lambda + \mu -1)/(\lambda \mu )\,,
$$
$$
|T_y|=\gamma \beta \left( (\lambda + \mu -1)/2 \right) \alpha (1- \alpha )\,.
$$
By using the program package GAP \cite{GAP}, 
$$
\begin{cases}
|(T^*)_y|= (1/2) \cdot (\lambda + \mu -1)^2
\cdot (1- \gamma ) \beta \alpha (1- \alpha )/ \\
\bigl[ \mu \cdot [ 1+(\lambda + \mu -1) \alpha ] \cdot 
[  \beta (1 -\alpha -\gamma
\alpha + \gamma  \alpha ^2) +  \\
\gamma \beta \alpha (1-\alpha ) \lambda + 
\alpha (1+ \gamma \beta -\gamma \alpha \beta ) \mu ] \bigr] 
\,.
\end{cases}
$$
Here the 
denominator is a product of three factors, all being positive. (For
the third factor observe that the coefficients of $ \lambda $, or $\mu $ are
non-negative or positive, respectively,
and the constant term is minimal for $ \gamma =1$,
and is then non-negative.) The second factor of the denominator will cancel
with $|C_i|$, 
and its third factor will be estimated from above as follows.
The coefficients of $\lambda $, or $ \mu $, in it, are 
estimated from above by setting $ \beta = \gamma =1$, and the constant term is
estimated from above by setting $\gamma =0$, and $\beta =1$. Thus we obtain
the upper estimates $ \alpha (1 - \alpha )$, or 
$ \alpha ( 2 - \alpha)$, or $1 - \alpha $, 
respectively. These can be further estimated from above by $1/4$, or $1$, or
$1$, respectively. 

Hence, minimizing for $\gamma  \in [0,1]$,
$$
\begin{cases}
|T_y| \cdot |(C_o)^*| +|C_i| \cdot |(T^*)_y| \ge \\
(1/4) \cdot \left( ( \lambda + \mu -1)^2 / \mu \right) \cdot
\min \, \{ 1/(1+\lambda /4 +\mu ), 1/ \lambda \}
\cdot \beta \alpha (1 - \alpha )
\end{cases}
$$
(the first term being estimated from below by setting $\alpha =1$ in the
denominator of the second factor of $|(C_o)^*|$).
Changing the roles of $\lambda, \mu $, 
of $\beta , 1- \beta $, and of $ \gamma ,
\delta $, we obtain similarly
$$
\begin{cases}
|T_z| \cdot |(C_o)^*| +|C_i| \cdot |(T^*)_z| \ge \\
(1/4) \cdot \left( ( \lambda + \mu -1)^2 / \lambda \right) \cdot 
\min \, \{ 1/(1+\lambda  +\mu /4), 1/ \mu \}
\cdot (1-\beta )\alpha (1 - \alpha )\,.
\end{cases}
$$
Hence, 
$$ 
\begin{cases}
|C| \cdot |C^*| \ge |C_i| \cdot |(C_o)^*| + 
(|T_y|+|T_z|) \cdot |(C_o)^*| +|C_i| \cdot (|(T^*)_y|+|(T^*)_z|) \\ 
\ge f( \lambda , \mu )+g( \lambda , \mu ) \cdot \alpha (1 - \alpha )\,.
\end{cases}
$$
$\blacksquare $


\begin{Cor}
\label{sectionstab3}
Under the hypotheses of Lemma~\ref{sectionstab2}, let 
$$
|C| \cdot |C^*| \le (1+\varepsilon ) \cdot f( \lambda , \mu )\,,
$$
where $\varepsilon \in \left( 
0, g(\lambda ,\mu ) / \left( 4f(\lambda , \mu ) \right) \right)
$. Further let
$\alpha _{ \pm } := \bigl[ 1 \pm  $ 
\newline
$ {\sqrt{1-\left( 4f(\lambda , \mu )/g(\lambda ,\mu )
\right) \varepsilon }} \bigr] /2$
and let \,$\alpha _+ + (1- \alpha _+) 
\min \, \{ (1 - \lambda )/ \mu, (1 - \mu ) / \lambda \} >0$.
Then 
$$
{\text{either\,\,\,\,}} 
C \subset [1+(\lambda + \mu -1) \alpha _-] \cdot [o,u,v], 
$$
$$
{\text{or\,\,\,\,}}
C \supset [\alpha _+ + (1- \alpha _+) \cdot 
\min \, \{ (1 - \lambda )/ \mu, (1 - \mu ) / \lambda \}] \cdot [o,u,v,p]\,.
$$ 
\end{Cor}


{\bf Proof.}
We use the notations from the proof of Lemma \ref{sectionstab2}.

By the hypotheses and Lemma~\ref{sectionstab2}, for $\alpha $ from 
Lemma~\ref{sectionstab2}, we have
$$ 
f( \lambda , \mu ) \cdot (1+ \varepsilon ) \ge |C| \cdot |C^*| \ge f( \lambda
, \mu )+g( \lambda , \mu ) \alpha (1- \alpha )\,,
$$
hence 
$$
\alpha ^2 - \alpha + \left( f( \lambda , \mu )/g( \lambda , \mu ) \right)
\varepsilon \ge 0\,,
$$
i.e., $\alpha \le \alpha _-$, or $\alpha \ge \alpha _+$, where $\alpha _{\pm}
\in (0, \infty )$ and
$\alpha _- < \alpha _+$.

Let $x \in C \cap [u,v,p]$, with $|[u,v,x]|$ maximal. 
Then $C$ lies below the line $l:=\{ y \mid y 
{\text{\rm{ lies above the line }}} uv, {\text{\rm{ and }}} |[u,v,y]|=
\alpha \cdot |[u,v,p]| \} $. 
If $\alpha \le \alpha _-$, then $C$ lies below the line $l_-$, defined
analogously to $l$, but using $\alpha _-$ rather than $\alpha $.
If $\alpha \ge \alpha _+$, then 
$C \supset [o,u,v,x]$, hence $C \supset [o,u,v, \nu x]$, where $ \nu x$ 
lies on the
line $l_+$, defined
analogously to $l$, but using $\alpha _+$ rather than $\alpha $.
Hence $C$ contains the quadrangle obtained from 
$[o,u,v, \nu x]$, 
by replacing its side lines $u( \nu x),v( \nu x)$ by lines through $ \nu x$,
parallel to $up,vp$. We further
diminish this last quadrangle by translating its side lines
parallel to $up$ or $vp$ so 
that they should contain the points of intersection
of the sides $vp$ or $up$ with the line $l_+$, respectively. 
The formulas in the corollary 
then follow by simple calculations.
$\blacksquare $

\hskip.0cm


{\bf{Remark.}}
It is probable that with more work one could sharpen the stability
estimates in the second proof of Lemma \ref{sectionstab},
iterating
further the construction of inscribed/circumscribed polygons (defining, in an
analogous manner, 
some closer approximations $C_i \subset C_i' \subset C_i''
\subset C \subset C_o'' \subset C_o' \subset C_o$, etc.).  
However, this way does not seem to be suitable to
give estimates which are sharp, up to a quantity $o( \varepsilon )$.

\hskip.0cm


The first inequality in the next lemma is related to
\cite{KR}, Proposition 1, but is formulated with
constants according to our particular needs in this paper. 
The second inequality in our next lemma is related to 
an opposite inequality as in Proposition 2 of \cite{KR}, but the idea of the
proof is similar.


\begin{Lem}
\label{centrestab}
Let $d \ge 2$ be an integer, $K_0 \subset {\mathbb R}^d$ be a convex body, 
and let
$0 < \varepsilon _1 \le \varepsilon _1(K_0):= \min \, \{ 1/2,
2^{-2d-1}
\left( \kappa _{d-1}/(d \kappa _d ^2) \right) \cdot |K_0|
/({\text{\rm{diam}}}\,K_0)^d \}
$. 
Let $K \subset {\mathbb R}^d$ be a convex body, and let 
$(1- \varepsilon _1) K_0 + a 
\subset K \subset (1+ \varepsilon _1)  K_0 +b$, where $a,b \in 
{\mathbb R}^d$. 
Then 
$$
\| s(K)-s(K_0) \| \le c_1(K_0) \cdot \varepsilon _1\,,
$$
where 
$$
c_1(K_0):=
({\text{\rm{diam}}}\,K_0)^{(d+1)^2} 
|K_0|^{-d-2} \cdot d (d \kappa _d / \kappa _{d-1}) ^{d+2} \,.
$$
If, moreover, 
$\varepsilon _2>0$ and $| K_0|
\cdot | \left( K_0-s(K_0) \right) ^*| \le | K|
\cdot \left( K-s(K) \right) ^*|$ and $c \in {\text{\rm{int}}}\,K$, and
$|K| \cdot |(K-c)^*| 
\le |K_0| \cdot
|\left( K_0-s(K_0) \right) ^*| + \varepsilon _2 \le \kappa _d ^2$, 
then 
$$
\|c-s(K_0)\|  \le  c_1(K_0) \cdot \varepsilon _1+ c_2(K_0) \cdot 
\sqrt{ \varepsilon _2 }\,,
$$
where
$$
c_2(K_0):=\sqrt{ ({\text{\rm{diam}}}\,K_0)^{d+2}/| K _0| } \cdot
\sqrt{2^{d+3} / \left( (d+1) \kappa _d \right) }\,.
$$
\end{Lem}


{\bf Proof.}
We will suppose that the point of homothety of $(1- \varepsilon _1) K_0 
+ a$ and $(1+ \varepsilon _1) K_0 +b$, that is in the first body, is
$o$ 
(this does not change $K_0-s(K_0)$, $K-s(K)$, $K-c$; namely, we consider $c$
as ``fixed to $K$'').
Thus $a=b=o$ can be supposed.

We have 
\begin{equation}
\label{derivatives}
\begin{cases}
|(K-x)^*|=(1/d) \int _{S^{d-1}} \left( h_K(u)- \langle u,x \rangle \right)
^{-d} du\,, \\
(\partial / \partial x_d) | (K-x)^* | = \int _{S^{d-1}} u_d \left(
h_K (u) - \langle u,x \rangle \right) ^{-d-1} du\,, \\
(\partial / \partial x_d)^2 | (K-x)^* | = (d+1) \int _{S^{d-1}} u_d^2 \left(
h_K (u) - \langle u,x \rangle \right) ^{-d-2} du \\ 
\ge (d+1) ({\text{diam}}\,K)^{-d-2} \kappa _d\,,
\end{cases}
\end{equation}
where $u=(u_1, \ldots , u_d)$, and  
$h_K$ is the support function of $K$, 
and $\kappa _d$ the volume of the unit ball in
${\mathbb{R}}^d$.
The analogues of these formulas hold for the first and second directional
derivatives in any direction.

First
we estimate $\delta := \| s(K)-s(K_0) \| $ 
from above. We may assume that $s(K)-s(K_0)
= (0, \ldots , 0, \delta )$, where $\delta > 0$. 

We begin by showing that
$s(K) \in {\text{int}} \left( (1- \varepsilon _1)K_0 \right) $, 
and even estimate ${\text{dist}} 
\left( s(K) , {\text{bd}}\,[(1- \varepsilon _1)K_0] \right) $ from below.
Let $ \eta := {\text{dist}} \left( s(K) , {\text{bd}}\,K  \right) \le $
${\text{dist}} 
\newline
\left( s(K) , {\text{bd}}\,[(1 + \varepsilon _1)K_0] \right) $. 
Then $\left( K-s(K) \right) ^* $ contains
$({\text{diam}}\,K)^{-1}B^d$, and a
point at distance $\eta ^{-1}$ from $o$ (with $B^d$ the unit ball about $o$).
Therefore 
\begin{equation}
\label{eta}
\kappa _d ^2 \ge | K | \cdot | \left( K-s(K) \right) ^* | 
\ge |K| \cdot ({\text{diam}}\,K)^{-d+1} ( \kappa _{d-1}/d ) \eta ^{-1}\,.
\end{equation}

Hence, by $\varepsilon _1 \le 1/2$, 
$$
\begin{cases}
\eta _0 := 2^{-2d+1}
\left( \kappa _{d-1}/(d \kappa _d ^2) \right) \cdot |K_0|/
[({\text{diam}}\,K_0) ^{d-1}] \le  \\
\left( \kappa _{d-1}/(d \kappa _d ^2) \right) \cdot |K|/
({\text{diam}}\,K) ^{d-1}
\le \eta \le  \\
{\text{dist}} \left( s(K) , {\text{bd}}\,[ (1 + \varepsilon _1)K_0 ] \right)
\le \\
{\text{dist}} \left( s(K) , {\text{bd}}\,[ (1 - \varepsilon _1)K_0 +2
\varepsilon _1 \cdot {\text{diam}}\,K_0 \cdot B^d ] \right) \,.
\end{cases}
$$
Thus, for $\varepsilon _1
\le \eta _0/(4 \cdot {\text{diam}}\,K_0)$, 
$$
s(K) \in {\text{int}}\,[(1 - \varepsilon _1)K_0] {\text{\,\,\,\,and\,\,\,\,}}
\eta _0/2 \le
{\text{dist}} \left( s(K) , {\text{bd}}\,[ (1 - \varepsilon _1)K_0 ] \right)
\,.
$$

Then, using convexity of the function $t^{-d-1}$ for $t>0$, and (\ref{eta})
for $K_0$, rather than $K$, we have
\begin{eqnarray}
\label{deltaepsilon}
\begin{cases}
0= \int _{S^{d-1}} u_d \left(
h_K (u) - \langle u, s(K) \rangle \right) ^{-d-1} du \ge \\
\int _{S^{d-1}} u_d [ 
h_{K_0} (u) + {\text{sg}} \,u_d \cdot
\varepsilon _1 h_{K_0} (u) - \langle u, s(K_0) \rangle 
- \delta u_d ] ^{-d-1} du \ge \\
\int _{S^{d-1}} u_d \left( h_{K_0}(u)- \langle u, s(K_0) \rangle \right) 
^{-d-1}du + \\
(d+1) \cdot \int _{S^{d-1}} u_d \left( \delta u_d- \varepsilon _1 \cdot 
{\text{sg}}\,u_d \cdot h_{K_0}(u) \right) \times \\ 
\left( h_{K_0}(u)- \langle u, s(K_0) \rangle \right) ^{-d-2}du \ge \\
\delta (d+1) ({\text{diam}}\,K_0)^{-d-2} \int _{S^{d-1}} u_d^2 du - \\
\varepsilon _1 (d+1)  \cdot {\text{diam}}\,K_0 \cdot [ \left(
\kappa _{d-1}/(d \kappa _d ^2) \right) | K_0 |/ \\
({\text{diam}}\,K_0)^{-d+1} ]^{-d-2} \int _{S^{d-1}} | u_d | du \,.
\end{cases}
\end{eqnarray}
Here, $\int _{S^{d-1}} u_d^2 du=\kappa _d$, and $\int _{S^{d-1}} | u_d | du
\le \int _ {S^{d-1}} du$, and comparing the first and last terms of 
(\ref{deltaepsilon}), we get the first inequality of the Lemma.

We turn to the second inequality. We have
\begin{equation}
\label{distcsantalo}
\|c-s(K_0)\| \le \| c-s(K) \| + \| s(K) - s(K_0) \| \le \| c-s(K) \| +
c_1(K_0) \varepsilon _1
\end{equation}
and
$$
|K| \cdot |(K-c)^*| \le 
|K_0| \cdot |\left( K_0-s(K_0) \right) ^*| + \varepsilon _2 \le
|K| \cdot |\left( K-s(K) \right) ^*| + \varepsilon _2 \,.
$$
We use (\ref{derivatives}) on the line $s(K)c$, 
which gives
\begin{equation}
\label{distcsantalok}
| K |(d+1) ({\text{diam}}\,K)^{-d-2} \kappa _d \cdot \| c-s(K) \|^2 /2 \le 
\varepsilon _2\,.
\end{equation}
(\ref{distcsantalo}) and (\ref{distcsantalok}) give the second inequality of
the Lemma.
$\blacksquare $

\vskip.3cm


{\bf Proof of Theorem~\ref{nfoldstab}.}
{\bf{1.}} First we estimate $\delta _{BM}^s (K,R_n)$ from above.
Here we may assume that $o$ is the Santal\'o point of $K$, i.e., its centre of
rotational symmetry.
As explained in \S 3, there exist
regular $n$-gons $K_i$ and $K_o$ centred at the origin,
such that $K_i\subset K \subset K_o$,
and the midpoints of the sides of $K_o$ are the vertices of $K_i$.
Assuming that the unit circular disc about $o$ is the incircle of
$K_o$, we have $K_o^*=K_i$. Now the radii from $o$ to 
the vertices of $K_i$
divide $K_o$ into $n$ congruent deltoids
$\widetilde{C}_1,\ldots,\widetilde{C}_n$ whose
common vertex is the origin. In particular,
$\widetilde{C}_j^*:=\widetilde{C}_j\cap K_i$
is the corresponding triangular sector of $K_i$, where
$j=1,\ldots,n$. For
the congruent sectors $C_j=\widetilde{C}_j\cap K$
of $K$, and the congruent sectors $C_j^*=\widetilde{C}_j\cap K^*$
of $K^*$, where $j=1,\ldots,n$, we have
$$
(1+\varepsilon)n^2|\widetilde{C}_1|\cdot |\widetilde{C}_1^*|
= (1+\varepsilon)|K_i|\cdot|K_o|\geq |K|\cdot|K^*|=
n^2|C_1|\cdot|C_1^*|.
$$
We observe that $\widetilde{C}_1^*=[o,u,v]$ and
$\widetilde{C}_1=[o,u,v,p]$, where
$p=\lambda u+\lambda v$ for $\lambda=[\cos(\pi/n)]^{-2}/2$,
and
$$
|C_1|\cdot|C_1^*|\leq (1+\varepsilon)|[o,u,v,p]|\cdot |[o,u^*,v^*]|.
$$
We deduce by Lemma~\ref{sectionstab}
that either $C_1\subset(1+\gamma\varepsilon)\widetilde{C}_1^*$, or
$(1+\gamma\varepsilon)^{-1}\widetilde{C}_1\subset C_1$, where
$\gamma:=6(1+\sqrt{2\lambda})\le 18 $.
Therefore the rotational symmetry yields that
either $K\subset (1+18 \varepsilon )K_i$, or
$(1+18 \varepsilon )^{-1}K_o\subset K$.

{\bf{2.}} Now we turn to the proof of the stability of the centre of polarity.
The point $x$ is the point of homothety of $\lambda _1 R_n+x$ and $\lambda _2
R_n+x$, and $x \in \lambda _1 R_n+x$. We will suppose $x=o$\,; then $o \in
R_n$. Simultaneously,
we have to replace $K^*$ with $(K-c)^*$, for some $c \in {\text{int}}\,K$
(``fixed to $K$'').
Let
$K_{0,n}:=
[(\lambda _1 +\lambda _2)/2]R_n$ (this will take over the role of $K_0$ from
Lemma \ref{centrestab}). Then $\lambda _1 R_n \subset K \subset
\lambda _2 R_n$ and $\lambda _2 / \lambda _1 \le 1+18 \varepsilon $ imply
\begin{equation}
\label{lambda12}
\begin{cases}
K_{0,n} (1-9 \varepsilon ) \subset K_{0,n}/[(1+\lambda _2 / \lambda _1)/2]
\subset K \\
\subset K_{0,n}/[(1+\lambda _1 / \lambda _2 )/2] \subset K_{0,n} (1+9
\varepsilon)\,.
\end{cases}
\end{equation}
Note that by hypothesis $\varepsilon < 1/18$, 
so here 
\begin{equation}
\label{epsilon}
1-9 \varepsilon > 1/2 \,\,\,\,(>0)\,.
\end{equation}

Restricting Lemma~\ref{centrestab} to $d=2$, we have
$\varepsilon _1(K_{0,n})=[1/(32 \pi ^2)] \cdot |K_{0,n}|/$
\newline
$({\text{diam}}\,K_{0,n})^2$
and $c_1(K_{0,n})=2 \pi ^4 \cdot ({\text{diam}}\,K_{0,n})^9|K_{0,n}|^{-4}$,
and $c_2(K_{0,n})=4 $
\newline
$\sqrt{2/(3 \pi )} \cdot ({\text{diam}}\,K_{0,n})^2|K_{0,n}|^{-1/2}$. 
Here ${\text{diam}}\,K_{0,n}=1$ (in the Euclidean norm mentioned 
in the Theorem), hence $\min _n |K_{0,n}|$ 
is attained for $n=3$. So
$\min _n  \varepsilon _1(K_{0,n})=\sqrt{3}/(128 \pi ^2)$, and
$\max _n c_1(K_{0,n})=512 \pi ^4 /9$, and $\max _n c_2(K_{0,n})$
\newline
$ =8\sqrt{2/\pi } 3^{-3/4}$.

We apply Lemma~\ref{centrestab} for $d=2$, replacing 
there
$\varepsilon _1(K_{0,n})$ by $\min _n  
\varepsilon _1(K_{0,n})$ and $c_1(K_{0,n})$ by $\max _n
c_1(K_{0,n})$, and $c_2(K_{0,n})$ by $\max _n 
c_2(K_{0,n})$. By 
(\ref{lambda12}) we may choose
$\varepsilon _1 := 9 \varepsilon $. Also, by the hypothesis of the theorem,
$\varepsilon _2$ can be chosen so that  
$\varepsilon _2 \le \pi ^2 \varepsilon  \,\,\left( \ge n^2 \sin ^2( \pi /n )
\right)$. So 
\begin{equation}
\label{cs(K_{0,n})1}
\|c-s(K_{0,n})\|  \le  (512 \pi ^4 /9) 
\cdot 9 \varepsilon + 
8\sqrt{2/\pi } 3^{-3/4} \cdot 
\pi \sqrt{\varepsilon }\,,
\end{equation}
for 
\begin{equation}
\label{cs(K_{0,n})2}
0 < \varepsilon \le \varepsilon ^* 
:=[\sqrt{3}/(128 \pi ^2)]/9 =0.0001523... < 1/18\,.
\end{equation}
However, we will use (\ref{cs(K_{0,n})1}) only for 
$0 < \varepsilon \le \varepsilon ^{**}$, for some
$\varepsilon ^{**} \in (0, \varepsilon ^*]$, 
to be chosen later.

First let $0 < \varepsilon \le \varepsilon ^{**}$.
Then (\ref{cs(K_{0,n})1}) gives
\begin{equation}
\label{cs(K_{0,n})3}
\|c-s(K_{0,n})\|  \le \left( 512 \pi ^4 
\sqrt{\varepsilon ^{**}}
+8\sqrt{2 \pi } 3^{-3/4} \right) \cdot
\sqrt{ \varepsilon }\,.
\end{equation}

Second let $\varepsilon \ge  \varepsilon ^{**}$.
Then we have $c \in $ int\,$K \subset \lambda _2 R_n +x= \lambda _2 R_n $, 
and $s(K_{0,n}) \in K_{0,n} =[(\lambda _1 + \lambda _2)/2]R_n \subset 
\lambda _2 R_n $ (the last inclusion following from $o \in R_n$). 
Hence, also using $\lambda _2 / \lambda _1 \le 2$, we have  
\begin{equation}
\label{cs(K_{0,n})4}
\begin{cases}
\|c-s(K_{0,n})\| \le {\text{diam}}\,(\lambda _2 R_n)=
\lambda _2 /[(\lambda _1+ \lambda _2)/2] \\
\le 4/3
\le [4/(3 \sqrt{\varepsilon ^{**}})] \cdot
\sqrt{ \varepsilon }\,.
\end{cases}
\end{equation}

By (\ref{cs(K_{0,n})3}) and (\ref{cs(K_{0,n})4}), 
we have
\begin{equation}
\label{cs(K_{0,n})5}
\|c-s(K_{0,n})\| \le \left( \max \, \{ 512 \pi ^4 
\sqrt{\varepsilon ^{**}}
+8\sqrt{2 \pi } 3^{-3/4}, 
4/(3 \sqrt{\varepsilon ^{**}}) \} \right)
\cdot \sqrt{\varepsilon }\,.
\end{equation}
Now we minimize the coefficient of
$\sqrt{\varepsilon }$ in (\ref{cs(K_{0,n})5}), 
that is a function of $\varepsilon ^{**}$.
This minimum occurs when the two terms 
under the maximum sign are equal, that 
occurs for $\varepsilon ^{**}=
0.0000258...\,$, and its value is $262.30682...\,$.
(Observe that in fact
$0 < \varepsilon ^{**} < \varepsilon ^*$).
$\blacksquare $

\vskip.3cm


For the proofs of Theorems \ref{MahlerReisnerstab} and \ref{MahlerMeyerstab},
we need a simple stability version of
the inequality between the arithmetic and geometric means.
If $n\geq 2$ and $0<a_1\leq\ldots\leq a_n$, then
\begin{eqnarray*}
\frac{a_1+\ldots+a_n}{n\cdot(a_1\cdot\ldots\cdot a_n)^{1/n}} &=&
\frac{(\sqrt{a_n}-\sqrt{a_1})^2+2\sqrt{a_1a_n}+\sum_{1<j<n}a_j}
{n\cdot(a_1\cdot\ldots\cdot a_n)^{1/n}}\\
&\geq &  
\frac{(\sqrt{a_n}-\sqrt{a_1})^2+n\cdot(a_1\cdot\ldots\cdot a_n)^{1/n}}
{n\cdot(a_1\cdot\ldots\cdot a_n)^{1/n}}\\
&\geq & 1+\frac1n\left(1-\sqrt{\frac{a_1}{a_n}}\right)^2.
\end{eqnarray*}
It follows that 
\begin{equation}
\begin{cases}
\label{arithgeostab}
\mbox{\ if \ } \varepsilon\geq 0 \mbox{ \ and \ }
(a_1+\ldots+a_n)/[n\cdot(a_1\cdot\ldots\cdot a_n)^{1/n})] \leq
1+\varepsilon, \\
\mbox{\ then \ }
a_j/a_k\geq 1-2\sqrt{n\varepsilon}
\mbox{ \ for any $1\leq j,k\leq n$}.
\end{cases}
\end{equation}

(It is easy to give the sharp version of this inequality. Fixing
$a_2+...+a_{n-1}$, the minimum occurs when $a_2=...=a_{n-1}$; let their common
value be $x \in [a_1,a_n]$. Then derivation w.r.t. $x$ gives for the minimum
that $x=(a_1+a_n)/2$. However, the formula given above will be more convenient 
to apply.)


{\bf Proof of Theorem~\ref{MahlerReisnerstab}.}
{\bf{1.}} First we estimate $\delta _{BM}(K,P)$ from above.
Here we may assume that $o$ is the Santal\'o point of $K$, i.e., its centre of
symmetry.
As explained in \S 3, after Theorem \ref{KiKo}, we may
assume that $K_i\subset K\subset K_o$,
where $K_o$ and $K_i$ are squares centered at $o$,
the midpoints of the sides of $K_o$ are the vertices of $K_i$, and the sides
of $K_o$ have length $2$.
In particular, $K_i$ and $K_o$ are polar to each other.
It also follows that $\delta_{BM}(K,P)\leq 2$,
and hence if $\varepsilon\geq 0.005$, then we are done.
Therefore we assume that  $\varepsilon< 0.005$.

Now, in a suitable coordinate system, 
$K_o$ can be dissected into four unit squares
$S^1_o:=[0,1] \times [0,1]$, $S^2_o:=[-1,0] \times [0,1]$, 
$-S^1_o$ and $-S^2_o$.
 We write $S^j_i=S^j_o\cap K_i$, and $C_j=S^j_o\cap K$
 and $C^*_j=S^j_o\cap K^*$
 for $j=1,2$, and hence 
Lemma~\ref{section} implies $|C_j|\cdot|C^*_j|\geq |S^1_i|\cdot|S^1_o|$
 for $j=1,2$. We deduce by the hypothesis 
 $|K|\cdot|K^*|\leq (1+\varepsilon)\cdot 8 $ and Lemma~\ref{section} that
\begin{eqnarray}
\label{MahReis}
\begin{cases}
 (1+\varepsilon) \cdot |S^1_i|\cdot|S^1_o|\geq 
 [(|C_1|+|C_2|)/2] \cdot [(|C^*_1|+|C^*_2|)/2] \\
\geq\sqrt{|C_1|\cdot|C_2|\cdot|C^*_1|\cdot|C^*_2|}\,, 
{\text{\,\,\,\,and\,\,\,\,}}
 |C_j|\cdot |C^*_j| \geq |S^1_i|\cdot|S^1_o|.
\end{cases}
\end{eqnarray}
In particular,
 $$
 |C_j|\cdot|C^*_j|\leq(1+\varepsilon)^2 \cdot |S^1_i|\cdot|S^1_o|
 \leq(1+2.005\varepsilon) \cdot |S^1_i|\cdot|S^1_o|
 \mbox{ \ for $j=1,2$}.
 $$

To apply Lemma~\ref{sectionstab}, 
we have $\lambda=\mu=1$ and $\gamma=6(1+\sqrt{2})<15$
 both in the cases of $C_1$ and $C_2$.
Therefore, for each of $j=1,2$,
either $C_j\subset (1+\gamma \cdot 2.005\varepsilon)S^j_i$,
or $(1+\gamma \cdot 2.005\varepsilon)^{-1}S^j_o\subset C_j$.
If both of $C_1$ and $C_2$ satisfy either the first 
or the second condition, then $\delta_{BM}(K,P)\leq 1+31\varepsilon$,
and we are done. Therefore
we suppose that $C_1\subset (1+\gamma \cdot 2.005\varepsilon)S^1_i$,
and $(1+\gamma \cdot 2.005\varepsilon)^{-1}S^2_o\subset C_2$,
and seek a contradiction.
We have $|C_1|\leq (1+\gamma \cdot 2.005\varepsilon)^2/2$,
and since the diagonal of $S^2_o$ not containing $o$
is a subset of $C_2$, we also have
$|C_2|\geq (1+\gamma \cdot 2.005\varepsilon)^{-1}$.
It follows by $\varepsilon<0.005$ that
$|C_1|<(1-2\sqrt{2\varepsilon})|C_2|$.
On the other hand, (\ref{arithgeostab}) applied in (\ref{MahReis}) leads to
$|C_1|\geq(1-2\sqrt{2\varepsilon})|C_2|$, a
contradiction. 

{\bf{2.}}
The stability of the centre of polarity is deduced from 
Lemma~\ref{centrestab}
like in Theorem~\ref{nfoldstab}, by supposing $x=0$.
Simultaneously,
we have to replace $K^*$ with $(K-c)^*$, for some $c \in {\text{int}}\,K$
(``fixed to $K$''). Let $K_0:=[(\lambda _1 + \lambda _2)/2]P$.
Now 
$\varepsilon _1 (K_0)=1/(64 \pi ^2)$,
and $c_1(K_0)=64 {\sqrt{2}}\pi ^4$, and
$c_2(K_0)=16\sqrt{2/(3 \pi )}$.
We only note that 
by hypothesis $\varepsilon < 0.005$, 
and then
we use the sharper estimate
$\delta _{BM}(K,P) = \delta _{BM}(K,K_0) \le 1 + 31 \varepsilon $.
Then, rather than (\ref{lambda12}) and (\ref{epsilon}), we have 
\begin{equation}
\label{lambda12'}
\begin{cases}
K_0 \left( 1-(31/2) \varepsilon \right) 
\subset K_0/[(1+\lambda _2 / \lambda _1)/2]
\subset K \\
\subset K_0/[(1+\lambda _1 / \lambda _2 )/2] \subset K_0 \left( 1+(31/2)
\varepsilon \right) \,,
\end{cases}
\end{equation}
for
\begin{equation}
\label{epsilon'}
1-(31/2) \varepsilon > 1-(31/2)/200 \,\,\,\,(>0)\,.
\end{equation}
By (\ref{lambda12'}) we may choose $\varepsilon _1:=(31/2) \varepsilon $, and
by hypothesis of the theorem we may choose $\varepsilon _2 :=8 \varepsilon $.

Then we have, analogously to (\ref{cs(K_{0,n})1}) and (\ref{cs(K_{0,n})2}), 
that 
\begin{equation}
\label{cs(K_{0,n})1'}
\|c-s(K_{0,n})\|  \le  (64 \sqrt{2} \pi ^4) 
\cdot (31/2) \varepsilon + 
16\sqrt{2/(3 \pi )} \cdot
{\sqrt{8}} \sqrt{\varepsilon }\,,
\end{equation}
for 
\begin{equation}
\label{cs(K_{0,n})2'}
0 < \varepsilon \le \varepsilon ^* 
:=[ \varepsilon _1 (K_0)]/(31/2) < 0.005 \,.
\end{equation}
Then, analogously to (\ref{cs(K_{0,n})1}), (\ref{cs(K_{0,n})4}) and 
(\ref{cs(K_{0,n})5}), we have, also using $\lambda _2 / \lambda _1 \le 2$, 
that for $0 < 
\varepsilon \le \varepsilon ^{**}$\,\,\,\,$(\le \varepsilon ^*)$ we have 
\begin{equation}
\label{cs(K_{0,n})3'}
\|c-s(K_0)\|  \le \left( 64 \sqrt{2} \pi ^4 (31/2) 
\sqrt{\varepsilon ^{**}}
+16\sqrt{2/(3 \pi )} \sqrt{8} \right) \cdot
\sqrt{ \varepsilon }\,,
\end{equation}
and for $\varepsilon \ge \varepsilon ^{**}$ we have
\begin{equation}
\label{cs(K_{0,n})4'}
\|c-s(K_0)\| 
\le [4/(3 \sqrt{\varepsilon ^{**}})] \cdot \sqrt{ \varepsilon }\,.
\end{equation}
Hence, for any $\varepsilon >0$, we have
\begin{equation}
\label{cs(K_{0,n})5'}
\begin{cases}
\|c-s(K_0)\| \le \\
\left( \max \, \{ 64 \sqrt{2} \pi ^4 (31/2) 
\sqrt{\varepsilon ^{**}}
+16\sqrt{2 /(3 \pi )} \sqrt{8},\,\,
4/(3 \sqrt{\varepsilon ^{**}}) \} \right)
\cdot \sqrt{\varepsilon }\,.
\end{cases}
\end{equation}
The optimal choice of $\varepsilon ^{**}$
is $\varepsilon ^*$. The distance to be estimated from above is at most
$335.10941... \cdot \sqrt{ \varepsilon }$.
$\blacksquare $

\vskip.3cm


{\bf Proof of Theorem~\ref{MahlerMeyerstab}.}
{\bf{1.}} First we estimate $\delta _{BM}(K,T)$ from above.
We may assume that $K$ is not a parallelogram,
and $o$ is the Santal\'o point of $K$.
As it is explained in \S 3, after Theorem \ref{KiKo}, we may
assume that $K_i\subset K\subset K_o$,
where $K_i$ and $K_o$ are regular triangles, and
the midpoints of the sides of $K_o$ are the vertices of $K_i$.
It also follows that $\delta_{BM}(K,T)\leq 4$,
and hence if $\varepsilon\geq 1/300$, then we are done.
Therefore we assume that  $\varepsilon< 1/300$.

We use the notation and ideas of the proof Theorem~\ref{KiKo};
in particular $o\in{\rm int}\,K_i$.
We may assume that the circumradius of $K_i$ is $1$,
and hence $d_1+d_2+d_3=3$, and $a=b=\sqrt{3}$.

Since $|K|\cdot|K^*|\leq (1+\varepsilon)\cdot |K_o|\cdot|K_o^*|$,
and we used the inequality
between arithmetic and geometric means
for $|C_1|,|C_2|,|C_3|$ in (\ref{sumCj}),
and for $d_1,d_2,d_3$
in the step from (\ref{sumproddj})
to (\ref{proddj}),  for $j,k=1,2,3$,
we deduce by (\ref{arithgeostab}) that
\begin{eqnarray}
\label{Cjstab}
|C_j|/|C_k|&\geq& 1-2\sqrt{3\varepsilon}\geq 4/5\,, {\rm{\,\,\,\,and}}
\\
\label{djstab}
d_j/d_k&\geq& 1-2\sqrt{3\varepsilon}\geq 4/5.
\end{eqnarray}
Since  $d_1+d_2+d_3=3$, we have
\begin{eqnarray}
\label{djlow}
d_j&\geq& \mbox{$3/(1+5/4+5/4)=6/7,$} {\rm{\,\,\,\,\,and}}
\\
\label{djup}
d_j&\leq& \mbox{$3/(1+4/5+4/5)=15/13$}.
\end{eqnarray}

Like in the proof of Theorem~\ref{MahlerReisnerstab}, by hypothesis, and
by Lemma~\ref{section},
$$
\begin{cases}
(1+\varepsilon)\left(\prod_{j=1}^3 ( |[o,x_{j-1},x_j,y_j]|\cdot
|[o,y_{j-1}^*,y_j^*]| ) \right)^{1/3} \ge \\
\left(\prod_{j=1}^3 ( |C_j|\cdot |C_j^*| ) \right)^{1/3}\,,
{\text{\,\,\,\,and}}
\end{cases}
$$
$$
|C_j|\cdot |C_j^*| \ge |[o,x_{j-1},x_j,y_j]|\cdot
|[o,y_{j-1}^*,y_j^*]|\,.
$$
Hence, for each $j=1,2,3$, we have
\begin{equation}
\label{3.1}
|C_j|\cdot |C_j^*|\leq (1+3.1\varepsilon)|[o,x_{j-1},x_j,y_j]|\cdot
|[o,y_{j-1}^*,y_j^*]|.
\end{equation}

Let $j=1,2,3$. To apply Lemma~\ref{sectionstab}, we define 
$\lambda_j,\mu_j>0$ by
$$
y_j=\lambda_jx_{j-1}+\mu_jx_j.
$$
Since $\lambda_j/\mu_j=|[o,x_j,y_j]|/|[o,x_{j-1},y_j]|=
d_j/d_{j-1}$,
(\ref{djstab}) implies
$$
\frac{\lambda_j+\mu_j}{\min \, \{\lambda_j,\mu_j\}}\leq \mbox{$1+5/4=9/4$}.
$$
Now the distances of $y_j$, or $o$ from the line through 
$x_{j-1},x_j$ are $3/2$,
or $\|x_{j+1}^*\|^{-1}=3/2-d_{j+1}\geq 9/26$, by (\ref{djup}), 
respectively, and hence
$$
\begin{cases}
\lambda_j+\mu_j=\langle x_j^*,y_j\rangle=
\langle x_j^*,x_j\rangle+
\| x_j^* \| \cdot \langle \| x_j^* \| ^{-1} x_j^*,y_j-x_j\rangle \\
\leq\mbox{$ 1+( 3/2 ) /( 9/26 ) =16/3$}.
\end{cases}
$$
We define $\gamma_j:=3[(\lambda_j+\mu_j)/\min \, \{\lambda_j,\mu_j\} ]
(1+\sqrt{\lambda_j+\mu _j})$,
and hence
$$
3.1\gamma_j\leq \mbox{$3.1\cdot 3\cdot (9/4) \cdot (1+4/\sqrt{3})<70. $}
$$
In particular, it follows by Lemma~\ref{sectionstab} and (\ref{3.1}) that
$$
\mbox{either\,\, $(1+70\varepsilon)^{-1}[o,x_{j-1},x_j,y_j]\subset C_j$\,,
\,\, or \,\,
$C_j\subset(1+70\varepsilon)[o,x_{j-1},x_j]$}\,.
$$

We note that $1+70\varepsilon\leq 5/4$ and $\|x_{j-1}-x_j\|=\sqrt{3}$.
If $(1+70\varepsilon)^{-1}[o,x_{j-1},$ $x_j,y_j]\subset C_j$, then
(\ref{djup}) yields
\begin{equation}
\label{Cjlarge}
\begin{cases}
|C_j|\geq |[o,x_{j-1},x_j,(4/5)y_j]|= \\
(4/5)\cdot (\sqrt{3}/2) \cdot (3/2-d_{j+1}+3/2)\geq 
(2\sqrt{3}/5)\cdot48/26>1.27.
\end{cases}
\end{equation}
On the other hand, if $C_j\subset(1+70\varepsilon)[o,x_{j-1},x_j]$,
then (\ref{djlow}) yields
\begin{equation}
\label{Cjsmall}
\begin{cases}
|C_j|\leq (5/4)^2 \cdot |[o,x_{j-1},x_j]|=
(5/4)^2 \cdot ({\sqrt{3}}/2) \cdot (3/2-d_{j+1})\leq \\
(5/4)^2 \cdot ({\sqrt{3}}/2)\cdot (9/14)<0.87.
\end{cases}
\end{equation}

Comparing (\ref{Cjstab}), (\ref{Cjlarge}) and (\ref{Cjsmall})
shows that either $(1+70\varepsilon)^{-1}[o,x_{j-1},x_j,y_j]$ 
\newline
$\subset C_j$
for all $j=1,2,3$, or
$C_j\subset(1+70\varepsilon)[o,x_{j-1},x_j]$ for all $j=1,2,3$.
Therefore either $(1+70\varepsilon)^{-1}K_o\subset K$,
or $K \subset (1+70\varepsilon)K_i$, and
hence the Banach-Mazur distance of $K$ from the triangles
is at most $1+70\varepsilon$.

{\bf{2.}}
The stability of the centre of polarity is deduced from 
Lemma~\ref{centrestab}
like in Theorem~\ref{nfoldstab} and 
Theorem~\ref{MahlerMeyerstab}, by supposing $x=o$. 
Simultaneously,
we have to replace $K^*$ with $(K-c)^*$, for some $c \in {\text{int}}\,K$
(``fixed to $K$''). Let $K_0:=[(\lambda _1 + \lambda _2)/2]T$.
Now 
$\varepsilon _1 (K_0)=\sqrt{3}/(128 \pi ^2)$ and $c_1(K_0)=512 \pi ^4 /9$
and $c_2(K_0)=8{\sqrt{2/(3 \pi )}} \cdot 3^{-1/4}$.
We only note that 
by hypothesis $\varepsilon < 1/300$, 
and then
we use the sharper estimate
$\delta _{BM}(K,T) = \delta _{BM}(K, K_0) \le 1 + 70 \varepsilon $.
Then, rather than (\ref{lambda12}), (\ref{lambda12'}) and 
(\ref{epsilon}), (\ref{epsilon'}), we have 
\begin{equation}
\label{lambda12''}
\begin{cases}
K_0 \left( 1-(70/2) \varepsilon \right) 
\subset K_0/[(1+\lambda _2 / \lambda _1)/2]
\subset K \\
\subset K_0/[(1+\lambda _1 / \lambda _2 )/2] \subset K_0 \left( 1+(70/2)
\varepsilon \right) \,,
\end{cases}
\end{equation}
for
\begin{equation}
\label{epsilon''}
1-(70/2) \varepsilon > 1-(70/2)/300 \,\,\,\,(>0)\,.
\end{equation}
By (\ref{lambda12''}) we may choose $\varepsilon _1:=(70/2) \varepsilon $, and
by hypothesis of the theorem, we may choose $\varepsilon _2 
:=(27/4) \varepsilon $.

Then we have, analogously to (\ref{cs(K_{0,n})1}), (\ref{cs(K_{0,n})1'}),
and (\ref{cs(K_{0,n})2}), (\ref{cs(K_{0,n})2'}) that 
\begin{equation}
\label{cs(K_{0,n})1''}
\|c-s(K_0)\|  \le  (512 \pi ^4 /9)
\cdot (70/2) \varepsilon + 
8\sqrt{2/(3 \pi )} \cdot 3^{-1/4}
{\sqrt{27/4}} \sqrt{\varepsilon }\,,
\end{equation}
for
\begin{equation}
\label{cs(K_{0,n})2''}
0 < \varepsilon \le \varepsilon ^* 
:=[ \varepsilon _1 (K_0)]/(70/2) < 1/300 \,.
\end{equation}
Then, analogously to (\ref{cs(K_{0,n})1}), (\ref{cs(K_{0,n})1'}), and
(\ref{cs(K_{0,n})4}), (\ref{cs(K_{0,n})4'}), and
(\ref{cs(K_{0,n})5}), (\ref{cs(K_{0,n})5'}) we have, 
also using $\lambda _2 / \lambda _1 \le 4$, 
that for $0 < 
\varepsilon \le \varepsilon ^{**}$\,\,\,\,$(\le \varepsilon ^*)$ there holds 
\begin{equation}
\label{cs(K_{0,n})3''}
\|c-s(K_0)\|  \le \left( ( 512 \pi ^4/9)(70/2) \sqrt{\varepsilon ^{**}}
+8\sqrt{2/(3 \pi )} \cdot 3^{-1/4} {\sqrt{27/4}} \right) \cdot
\sqrt{ \varepsilon }\,,
\end{equation}
and for $\varepsilon \ge \varepsilon ^{**}$ we have

\begin{equation}
\label{cs(K_{0,n})4''}
\|c-s(K_0)\| \le 
\lambda _2 /[(\lambda _1+ \lambda _2)/2] 
\le [8/(5 \sqrt{\varepsilon ^{**}})] \cdot \sqrt{ \varepsilon }\,.
\end{equation}

Hence, for any $\varepsilon >0$, we have
\begin{equation}
\label{cs(K_{0,n})5''}
\begin{cases}
\|c-s(K_0)\| \le \\
\left( \max \, \{ (512 \pi ^4/9)(70/2)
\sqrt{\varepsilon ^{**}} + 
8\sqrt{2 /(3 \pi )} 3^{-1/4} \sqrt{27/4},\,\,
8/(5 \sqrt{\varepsilon ^{**}}) \} \right) \\
\cdot \sqrt{\varepsilon }\,.
\end{cases}
\end{equation}
The optimal choice of $\varepsilon ^{**}$
is $\varepsilon ^*$. 
The distance to be estimated from above is at most
$916.69531... \cdot \sqrt{ \varepsilon }$.
$\blacksquare $


\vskip.3cm

We turn to the proof of Theorem~\ref{Egglestonstab}. We proceed analogously as
in Lemma~\ref{sectionstab2} and Corollary~\ref{sectionstab3}.
Again, the proof of Lemma~\ref{Egglestonlemma}
will use an idea of Behrend, \cite{Be}, proof of (77),
pp. 739-740, and of (112), pp. 746-747.

As in the proof of Theorem~\ref{MahlerMeyerstab}, we assume that $K_i \subset
K \subset K_o$, where $K_o=[a,b,c]\,$ and 
$K_i=[a',b',c']$ are regular triangles, and
$a'=(b+c)/2\,$, $b'=(c+a)/2\,$, $c'=(a+b)/2$. Now we assume $\| a-b \| =2$.
We let $\alpha _1:= \max \, \{ | [x,b',c']|/|[a,b',c'] | \mid x \in K \cap
[a,b',c'] \} $, and
$ \alpha _2:= \max \, \{ | [x,c',a']|/$
\newline
$|[b,c',a'] | \mid x \in K \cap
[b,c',a'] \} $, and
$\alpha _3:= \max \, \{ | [x,a',b']|/|[c,a',b'] | \mid x \in K \cap
[c,a',b'] \} $. 
Then $ \alpha _i \in [0,1]$, and we let $\alpha :=(\alpha _1+ \alpha _2+ \alpha
_3)/3 \in [0,1]$.


\begin{Lem}
\label{Egglestonlemma}
With the above notations, we have
$$
|K| \cdot | [ (K-K)/2 ] ^* | \ge 6 + (3/2) \alpha (1- \alpha ) \,.
$$
\end{Lem}

{\bf Proof.}
The supporting lines of
$K$, parallel to and different from the side lines of $K_o$, contain points
$a'',b'',c''$ of $K$, with $a''$ lying in the triangle $b'ac'$, etc.
We let $K_i':=[a',c'',b',a'',c',b'']$, and let
$K_o'$ be the hexagon bounded by all supporting lines of $K$ parallel to the
sides of $K_o$.
We have 
$$
K_i' \subset K \subset K_o'\,.
$$
Hence, 
\begin{equation}
\label{Eggleston1}
|K| \cdot | [(K-K)/2]^* | \ge | K_i'| \cdot | [(K_o-K_o)/2)]^* | \,. 
\end{equation}
Here 
\begin{equation}
\label{Eggleston2}
| K_i'|=({\sqrt{3}}/4) (1+ \alpha _1+ \alpha _2+ \alpha _3)\,, 
\end{equation}
\begin{equation}
\label{Eggleston3}
\begin{cases}
| [(K_o'-K_o')/2]^* | = 
2 ( 4/ \sqrt{3} ) ^2 [ (1+ \alpha _1)^{-1}(1+ \alpha _2)^{-1} + \\
(1+ \alpha _2)^{-1}(1+ \alpha _3)^{-1}+(1+\alpha _3)^{-1}(1+ \alpha _1)^{-1} 
] \sin (\pi /3)/2\,. 
\end{cases}
\end{equation}

Now, (\ref{Eggleston2}), (\ref{Eggleston3}), and 
the arithmetic-geometric mean inequality imply 
\begin{equation}
\label{Eggleston4}
\begin{cases}
| K_i'| \cdot | [(K_o-K_o)/2]^* | = 
2(1+3 \alpha )(3+3 \alpha ) \times \\
(1+ \alpha _1)^{-1}(1+\alpha _2)^{-1}(1+\alpha _3)^{-1} \ge 
6(1+3 \alpha )(1+ \alpha )^{-2} 
\,.
\end{cases} 
\end{equation}
Taking in account (\ref{Eggleston1}),
it suffices to show that the last quantity in (\ref{Eggleston4})
is at least 
$$
6+ (3/2) \alpha (1- \alpha )\,.
$$
However, if we replace in the last expression 
$3/2$ by some $c \ge 0$, this claimed inequality
becomes equivalent to
$$
\alpha (1-\alpha ) \left( 1-(c/6) (1+\alpha )^2 \right) \ge 0\,,
$$
that is (just) satisfied for $c=3/2$.
$ \blacksquare $


\vskip.3cm

{\bf Proof of Theorem~\ref{Egglestonstab}.}
We will use the notations in Lemma~\ref{Egglestonlemma} and its proof.
By hypotheses and Lemma~\ref{Egglestonlemma}, 
$$ 
6 \cdot (1+ \varepsilon ) \ge |K| \cdot | [(K-K)/2]^* | \ge 6
+(3/2) \alpha (1- \alpha )\,,
$$
hence 
$$
\alpha ^2 - \alpha + 4 \varepsilon \ge 0\,,
$$
i.e., $\alpha \le \alpha _-$, or $\alpha \ge \alpha _+$,
where $\alpha _{\pm} $
are the roots of the last polynomial. They 
are real, with $\alpha _- < \alpha
_+$, for 
$$
\varepsilon \in [0, 1/16)\,,
$$ 
which last inequality will be supposed preliminarily.

For $\alpha \le \alpha _-$ we have 
\begin{equation}
\label{alpha _-1}
\delta _{BM}(K,T) \le 1+\alpha _1+ \alpha _2+ \alpha _3 = 1 + 3 \alpha _- \,.
\end{equation}

Now let $\alpha \ge \alpha _+$. We proceed analogously, as in the proof of
Corollary~\ref{sectionstab3}.
We write $\beta _i:=1 - \alpha _i \in [0,1]$,
and $\beta :=1 - \alpha \in [0,1]$. Then $\beta = (\sum \beta _i)/3
\le \alpha _-$, hence $\beta _i \le \sum \beta _i \le 3 \alpha _-$. We have 
$K \cap [a',b',c] \supset [a',b',c'']$. We diminish this last triangle by 
retaining its side line $a'b'$, and replacing its sides
$[a',c''],[b',c'']$ by sides containing $c''$, and parallel to
$[a',c],[b',c]$. We further diminish this last triangle by retaining its side
line $a'b'$, and translating its side lines parallel to $[a',c],[b',c]$,
so that they should contain $c''(b'),c''(a')$, where $[c''(a'),c''(b')] \ni
c''$ is a chord of $[a',b',c]$ parallel to $[a',b']$, 
with $c''(a') \in [c,a']$ and 
$c''(b') \in [c,b']$. Of course this is only possible for $\beta _3 \le 
3 \alpha _- \le 1/2$;
therefore we preliminarily suppose also 
$$
\alpha _- \le 1/6\,, {\text{\,\,\,\,or, equivalently,\,\,\,\,}}
\varepsilon \in [0, 1/28.8] \,\,\left( \subset [0,1/16) \right) \,.
$$
In this case the vertex $c'''$ 
of the last triangle opposite its side on $[a',b']$ depends only on $\beta _3$:
it lies on the angle bisector of the triangle $[a',c,b']$ at $c$, and 
$\| c'''-c \|= \beta _3 \sqrt{3}$. Lastly
we replace $c'''$ by $c''''$, which is constructed analogously as $c'''$, but
replacing at the beginning $\beta _3 $ by $3 \alpha _- \,\,( \ge \beta _3 )$.
Analogously we define the points $a'''',b''''$. 
Then $[a'''',b'''',c''''] \subset
[a',c'''',b',a'''',c',b''''] \subset K$, hence 
\begin{equation}
\label{alpha _-2}
\delta _{BM} (K,T) \le 1/\left( 1-(9/2) \alpha _- \right) \,.
\end{equation} 
Here we have $1-(9/2) \alpha _- \ge 1/4$, i.e., $\alpha _- \le 1/6$,
thus $a'''',b'''',c'''' \not\in {\text{int}}\,[a',b',$ 
\newline
$c']$.

Now, (\ref{alpha _-1}) and (\ref{alpha _-2}) give
\begin{equation}
\label{alpha _-3}
\begin{cases}
\delta _{BM} (K,T) \le  \max \, 
\{ 1+3 \alpha _-, \,\, 1/(1-(9/2) \alpha _-) \} =
\\
1/(1-(9/2) \alpha _-) = 
1+[(9/2) \alpha _- ]/[1-(9/2) \alpha _-]  \le \\
1+[(9/2) \alpha _- ]/[1-(9/2) (1/6)] =1+18 \alpha _- \,.
\end{cases}
\end{equation}
By convexity of the respective function,
\begin{equation}
\label{alpha _-4}
\alpha _- = (1-\sqrt{1-16 \varepsilon})/2 \le (24/5) \varepsilon , 
{\text{\,\,\,\,for\,\,\,\,}} \varepsilon \in [0, 1/28.8]\,.
\end{equation}
Thus, by (\ref{alpha _-3}) and (\ref{alpha _-4}),
\begin{equation}
\label{alpha _-5}
\delta _{BM}(K,T) \le 1+18 \alpha _- \le 1+86.4 \varepsilon \,.
\end{equation}
 
There remains the case $\varepsilon \ge 1/28.8$.
Then
\begin{equation}
\label{alpha _-6}
\delta _{BM}(K,T) \le 4 \le 1+86.4 \varepsilon \,.
\end{equation}
Lastly, (\ref{alpha _-5}) and (\ref{alpha _-6}) together prove the theorem.
$ \blacksquare $


\section{A short proof of the inequality of Mahler-Reisner} 

\begin{Thm}
\label{MahlerReisner} (Mahler-Reisner {\rm{\cite{Mah38}, \cite{R86}}}).
If $K$ is an $o$-symmetric convex body in $\R^2$, then
$$
|K|\cdot|K^*| \ge 8\,,
$$
with equality if and only if $K$ is a parallelogram.
\end{Thm}


In the proof of this theorem we use the results of \cite{MR06}, 
more exactly, the proof of
their Theorem 15. Actually we will make only slight modifications in its
proof.

\vskip.3cm


{\bf Proof of Theorem~\ref{MahlerReisner}.}
Like in \cite{MR06}, proof of their Theorem 15,
we may suppose that a diameter of $K$ 
(a segment of maximal length contained in $K$)
is $[(-1 ,0),\,\,(1 ,0)]$, 
where $K
\subset {\mathbb R}^2$ has a minimal volume product among $0$-symmetric convex
bodies. Let 
$$
K= \{ (x,y) \mid x \in [-1 ,\,\, 1 ],\,\,-f(-x) \le y \le f(x) \} \,,  
$$
where 
\begin{equation}
\label{MR}
\begin{cases}
f(x) {\text{\,\,\,\,is a concave function on\,\,\,\,}} [-1,\,\,1], 
{\text{\,\,\,\,with}} \\
f(-1)=f(1)=0, {\text{\,\,\,\,and\,\,\,\,}} f(x)>0 
{\text{\,\,\,\,for\,\,\,\,}} x \in (-1,\,\,1)\,. 
\end{cases}
\end{equation}
If the graph of $f$ consists of two segments, we
are done. If not, then, by Lemma 14 of \cite{MR06}, there are functions $g,h$,
both satisfying (\ref{MR})
above, both not proportional to $f$, such that $f=(g+h)/2$. 

Let $t \in [-1,\,\,1]$. 
Let $f_t:=f+t(h-g)/2$. Then the area of the convex body 
$$
K_t:=
\{ (x,y) \mid x \in [-1 , \,\,1 ],\,\,-f_t(-x) \le y \le f_t(x) \} 
$$
is a linear function of $t$.
By Theorem 1 of \cite{MR06} the reciprocal of 
$\varphi (t):=| [K_t-s(K_t)]^* | $ is a convex function of
$t$. Hence $\min \varphi $ is attained either for $t=-1$ or for $t=1$. Since
$K=K_0$ has minimum volume product, $\varphi $ is constant. Then, by
Proposition 7 of \cite{MR06}, $K_1$
is an affine image of $K$, by an affinity of the form $(x,y) \to (x, 
ux+vy+w)$. By \cite{MR06}, p. 140, Remark to Proposition 7, we have 
$h(x)=vf(x)+ux+w$. Putting here $x = \pm 1$, we see $u=w=0$. Hence, $h$ is
proportional to $f$, a contradiction. 
$ \blacksquare $

\vskip.3cm

{\bf{Acknowledgement.}} We express our thanks to the referee, whose careful
reading of our paper, and remarks about our paper 
substantially improved the presentation of our paper.
We also express our thanks to M. Kozachok and A. Magazinov for pointing out to
us the paper of Hanner.


\end{document}